


\input amstex
\documentstyle{amsppt}
\nologo
\pageheight{7.6in}
\pagewidth{5.6in}
\magnification=\magstephalf
\NoBlackBoxes
\topmatter
\pageno=-1
\leftheadtext{A conformally invariant sphere theorem in four dimensions}
\rightheadtext{S.-Y. A. Chang, M. J. Gursky, and P. C. Yang}
\TagsOnRight
\title{A conformally invariant sphere theorem
\\ in four dimensions}\endtitle

\vglue -6pc
\author{Sun--Yung A. Chang$^{\ssize 1}$, Matthew J. Gursky$^{\ssize 2}$,\\ 
and Paul C. Yang$^{\ssize 3}$ \\ \\ }\endauthor
\adjustfootnotemark{1}\footnotetext{Princeton University, Department of Mathematics, Princeton, NJ 08544-1000, supported by NSF Grant DMS--0070542.}
\adjustfootnotemark{1}\footnotetext{University of Notre Dame, 
Department of Mathematics, Notre Dame, IN  46556, supported in part by NSF Grant DMS--9801046 and 
an Alfred P. Sloan Foundation Research Fellowship.}
\adjustfootnotemark{1}\footnotetext{Princeton University, Department of Mathematics, Princeton, NJ  08544-1000, supported by NSF Grant DMS--0070526.}
\openup1.5\baselineskip
\toc \nofrills{Table of Contents}
\widestnumber\head{10.}
\head 0.  Introduction$\dotfill 1$\endhead
\head  1.  An existence result$\dotfill 6$\endhead
\head  2.  The proof of Theorem A$\dotfill 12$\endhead
\head  3.  A Weitzenb\"ock formula for Bach-flat metrics$\dotfill 12$\endhead
\head  4.  The proof of Theorem C$\dotfill 17$\endhead
\head \hskip .15cm  Appendix $\dotfill 24$ \endhead
\head \hskip .4cm References$\dotfill 33$\endhead
\endtoc
\vskip .2in
\endtopmatter

\parskip=5.5pt
\parindent=.7cm

\abstract{TO COME}\endabstract
\voffset=.40in

\def\bcw{\mathbin{\bigcirc\mkern-15mu\wedge}}

\define\vba{|\vec{B}|^2}
\define\vxa{|\vec{X}|^2}
\define\vya{|\vec{Y}|^2}
\define\si3{\sum_{i=1}^3}

\define\trdo{\nabla_0}

\define\bci2{\bigcirc  \hskip -.3cm \land }

\define\Rn{{\bold R}^n}
\define\R{\bold R}
\define\s2a{\sigma_2(A)}
\def\intl{\int \kern-8.9pt\hbox{$\diagup$}}
\def\intll{\int\kern-9.9pt\hbox{$\diagup$}}
\def\la123{\lambda_1^\pm \lambda_2^\pm \lambda_3^\pm}
\def\b12b22b32{b_1^2 \leq b_2^2 \leq b_3^2}
\def\bB1bB2bB3{b_1 \leq b_2 \leq b_3}
\vskip 1.5cm
\NoBlackBoxes

\document

\baselineskip=13pt

\newpage
\pageno=1
\vskip .3cm

\head{\bf 0.   Introduction}\endhead

Under what conditions on the curvature can we conclude that a smooth, closed 
Riemannian manifold is diffeomorphic (or homeomorphic) to the sphere?
A result which addresses this question is usually referred to as a 
{\it sphere theorem}, and the literature abounds with examples (see 
Chapter 11 of [Pe] for a brief survey).

In this paper we concentrate on four dimensions, where Freedman's 
work is obviously very influential.  For example, any curvature 
condition which implies the vanishing of the de Rham cohomology 
groups $H^1(M^4,\bold R )$ and $H^2(M^4,\bold R )$ will, by Freedman's
result ( \cite {Fr}), imply that $M^4$ is covered by a homeomorphism sphere.  
At the same time, there are some very interesting results
which characterize the {\it smooth} four-sphere.  

An example of particular importance to us is the work of Margerin 
(\cite {Ma2}), in which he formulated a notion of ``weak curvature pinching.''
To explain this we will need to establish some notation. 
Given a Riemannian four-manifold $(M^4,g)$, let $Riem$ denote the curvature
tensor, $W$ the Weyl curvature tensor, $Ric$ the Ricci tensor, 
and $R$ the scalar curvature.  The usual decomposition of $Riem$ under
the action of $O(4)$ can be written   
$$
Riem = W + {\tsize{1\over 2}} E  \bcw g + {\tsize{1\over 24}} R g \bcw g,
\tag0.1
$$
where $E = Ric - {1 \over 4} Rg$ is the trace-free Ricci tensor and 
$\bcw$ denotes the Kulkarni-Nomizu product.  If we let 
$Z = W + {1 \over 2} E \bcw g$, then 
$$
Riem = Z + {\tsize{1 \over 24}} R g \bcw g.
$$
Note that $(M^4,g)$ has constant curvature if, and only if, $Z \equiv 0$.  
We now define the scale-invariant "weak pinching" quantity
$$
WP \equiv {|Z|^2 \over R^2} = {|W|^2 + 2|E|^2 \over R^2} 
\tag0.2
$$
where $|Z|^2 = Z_{ijkl}Z^{ijkl}$ denotes the norm of $Z$ 
viewed as a $(0,4)$--tensor.   

Margerin's main result states that if $R > 0$ and $WP < {1 \over 6}$, 
then $M^4$ is diffeomorphic to either $S^4$
or ${\bold R}P^4$.
Moreover, this "weak pinching" condition is {\it sharp}:  
The spaces $({\bold C}P^2, g_{FS})$ and 
$(S^3 \times S^1, g_{prod.})$ both have $R > 0$ and $WP \equiv {1 \over 6}$. 
Indeed, using a holonomy reduction argument, Margerin also proved the
converse, in the sense that these manifolds (and any quotients) 
are characterized by the property that $WP \equiv {1 \over 6}$.  

Margerin's proof relied on an important tool in the subject 
of sphere theorems, namely, Hamilton's Ricci flow.  In fact, previously
Huisken (\cite {Hu}) and Margerin (\cite {Ma1}) independently
had used the Ricci flow to 
prove a similar pinching result, but with a slightly worse constant.
In addition, Hamilton (\cite {Ha}) had used his flow to study 
four-manifolds with positive curvature operator.  As Margerin
points out in his introduction, there is no relation between weak 
pinching and positivity of the curvature operator; indeed,
weak pinching even allows for some negative sectional curvature.

On the homeomorphism level, both Margerin's and Hamilton's 
curvature assumptions already imply that the underlying manifold
is covered by a homeomorphism sphere.  If $(M^4,g)$ has positive 
curvature operator, then the classical Bochner theorem
implies that $H^1(M^4,\bold R) = H^2(M^4,
\bold R) = 0$. 
As we observed above, Freedman's work then gives the 
homeomorphism type of the cover of $M^4$.  For manifolds satisfying 
$WP < {1 \over 6}$ the vanishing of harmonic forms
is less obvious, but does follow from \cite {Gu}. 
 
One drawback to the sphere theorems described above is that 
they require one to verify a {\it pointwise} condition on the curvature. 
In contrast, consider the (admittedly much simpler) case of surfaces. 
For example, if the Gauss curvature of the surface $(M^2,g)$ satisfies 
$\int K dA > 0$, then $M^2$ is diffeomorphic
to $S^2$ or $\bold R P^2$. 
In addition to this topological classification, the uniformization 
theorem implies that $(M^2,g)$ is conformal
to a surface of constant curvature, which is then covered isometrically by
$S^2$.  Therefore, in two dimensions one has a "sphere theorem" which only
requires one to check an {\it integral} condition on the curvature.  

Our goal in this paper is to generalize this situation, by
showing that 
the smooth four-sphere is also characterized by an 
integral curvature condition.  As we shall 
see, our condition has the additional properties of being {\it sharp} and
{\it conformally invariant}.  Although there are different---though 
equivalent---ways of stating our main result, the simplest version involves
the Weyl curvature and the Yamabe invariant:

\proclaim{Theorem A} 
Let $(M^4,g)$ be a smooth, closed four-manifold for which 
\roster 
{\item"{(i)}" the Yamabe invariant $Y(M^4,g) > 0$, and
\smallskip
\item"{(ii)}" the Weyl curvature satisfies
$$
\int_{M^4} |W|^2 dvol < 16 \pi^2 \chi(M^4).
\tag0.3
$$
\endroster
\smallskip\noindent
Then $M^4$ is diffeomorphic to either $S^4$
or $\bold R P^4$.
\endproclaim

\remark{Remarks}
\roster
\item"{1.}"  Recall that the Yamabe invariant is defined by
$$
Y(M^4,g) = \inf_{\tilde g \in [g]} 
    vol(\tilde g)^{-{1\over 2}}\int_{M^4} R_{\tilde g} dvol_{\tilde g},
$$
where $[g]$ denotes the conformal class of $g$. 
Positivity of the Yamabe invariant implies that
$g$ is conformal to a metric of strictly positive scalar curvature.  

\item"{2.}"  In the statement of Theorem A, the norm of the Weyl 
tensor is given by $|W|^2 = W_{ijkl}W^{ijkl}$;
i.e., the usual definition when $W$ is viewed as a section of 
$\otimes^4 T^{*}M^4$.  However, if one views $W$ as a section of 
$End(\Lambda^2(M^4))$, then the convention is 
$|W|^2 = {1\over4}W_{ijkl}W^{ijkl}$.
This can obviously lead to confusion not only when comparing 
formulas from different sources, but also
when the Weyl tensor is interpreted in different ways in the 
same paper (which will be the case here).  To avoid this problem,
our convention will be to denote the $(0,4)$-norm using 
$| \cdot |$, and the $End(\Lambda^2(M^4))$-norm
by $\| \cdot \|$, which has the added advantage of emphasizing 
how we are viewing the tensor in question.
We should note that some authors (for example, Margerin) 
avoid this confusion by just defining an isomorphism between 
$\otimes^4 T^{*}M^4$ and $End(\Lambda^2(M^4))$ which induces 
the same norm on both.  But in our case, we will
adopt the usual identification: if $Z \in \otimes^4 T^{*}M^4$, 
then we identify this with $Z \in End(\Lambda^2(M^4))$ by 
defining $Z(\omega)_{ij} = {1\over 2}Z_{ij}^{\quad kl}\omega_{kl}$.

\item"{3.}"  That Theorem A relies on an integral curvature 
condition indicates the possibility one
could attempt to formulate a version which imposed much weaker 
regularity assumptions on the metric.
For example, in \cite {CMS}, Cheeger, M\"uller, and Schrader 
defined a notion of curvature on piecewise flat
spaces, which also allowed a generalization of the 
Chern--Gauss--Bonnet formula.
\endroster

By appealing to the Chern--Gauss--Bonnet formula, 
it is possible to replace (0.3) with a 
condition which does not involve the Euler characteristic.  Since
$$
8 \pi^2 \chi(M^4) = \int_{M^4} \bigl( {\tsize{1\over 4}}|W|^2 - 
    {\tsize{1\over 2}} |E|^2 + {\tsize{1\over 24}}R^2 \bigr)dvol,   
\tag0.4
$$
Theorem A is equivalent to

\proclaim{Theorem A$^{\prime}$} 
Let $(M^4,g)$ be a smooth, closed manifold for which 
\roster 
\item"{(i)}" the Yamabe invariant $Y(M^4,g) > 0$, and
\smallskip
\item"{(ii)}" the curvature satisfies
$$
\int_{M^4} \bigl( -{\tsize{1\over2}} |E|^2 
      + {\tsize{1\over 24}}R^2 - {\tsize{1\over4}} |W|^2 \bigr) dvol > 0.
\tag0.5
$$
\endroster
\smallskip\noindent
Then $M^4$ is diffeomorphic to either $S^4$
or $\bold R P^4$.
\endproclaim

Formulating the result of Theorem A in this manner allows us to 
explain the connection with the work of Margerin. 
This connection relies on recent work (\cite {CGY1},\cite{CGY2})
in which we established the existence of solutions to 
a certain fully nonlinear equation in conformal geometry. 
The relevance of this PDE work
to the problem at hand is explained in Section 1.    
Simply put, the results of \cite {CGY1} and \cite {CGY2} allow us 
to prove that under the hypotheses of Theorem A$^{\prime}$, 
there is a conformal metric for which the integrand in (0.5) is
{\it pointwise} positive.  That is, through a conformal deformation of metric,
we are able to pass from positivity in an average sense to pointwise positivity.
Now, any metric for which the integrand in (0.5) is positive must satisfy
$$
|W|^2 + 2|E|^2 < {\tsize{1\over 6}}R^2,
$$
by just rearranging terms.
Note that this implies in particular that $R > 0$. 
Dividing by $R^2$, we conclude that
$$
WP = {|W|^2 + 2|E|^2 \over R^2} < {1\over6}.
$$
The conclusion of the theorem thus follows from Margerin's work.

As we mentioned above, Theorem A is {\it sharp}. 
By this we mean that we can
precisely characterize the case of equality:

\proclaim{Theorem B} 
Let $(M^4,g)$ be a smooth, closed manifold which is not 
diffeomorphic to either $S^4$ or $\bold R P^4$.  Assume in addition that 
\roster 
\item"{(i)}" the Yamabe invariant $Y(M^4,g) > 0$,
\smallskip
\item"{(ii)}" the Weyl curvature satisfies
$$
\int_{M^4} |W|^2 dvol = 16 \pi^2 \chi(M^4).
\tag0.6
$$
\endroster
\smallskip\noindent
Then one of the following must be true:
\roster
\item"{1}" $(M^4,g)$ is conformal to $\bold C P^2$ with the 
Fubini-Study metric $g_{FS}$, or
\smallskip
\item"{2}" $(M^4,g)$ is conformal to a manifold which 
is isometrically covered by $S^3 \times S^1$
endowed with the product metric $g_{prod.}$.
\endroster

\endproclaim

The proof of Theorem B relies on a kind of vanishing result, 
in a sense which we now explain.
Suppose $(M^4,g)$ satisfies that hypotheses of Theorem B. 
If there is another metric in a small neighborhood of 
$g$ for which the $L^2$-norm of the Weyl tensor is smaller; i.e., 
$$
\int_{M^4} |W|^2 dvol < 16 \pi^2 \chi(M^4),
$$
then by Theorem A we would conclude that $M^4$ is diffeomorphic 
to either $S^4$ or $\bold R P^4$.
This, however, contradicts 
one of the assumptions of Theorem B.  Therefore, for every metric in
some neighborhood of $g$, 
$$
\int_{M^4} |W|^2 dvol \geq 16 \pi^2 \chi(M^4).
$$
Consequently, $g$ is a critical point (actually, a local minimum) 
of the {\it Weyl functional} $g \mapsto \int |W|^2 dvol$. 
The gradient of this functional is called the {\it Bach
tensor}, and we will say that critical metrics are 
{\it Bach--flat}.  Note that the conformal invariance of 
the Weyl tensor implies that Bach-flatness is a conformally
invariant property.  In fact, the Bach tensor is conformally 
invariant, \cite {De}.

Theorem B is then a corollary of the following classification 
of Bach-flat metrics:

\proclaim{Theorem C} 
Let $(M^4,g)$ be a smooth, closed manifold which is not 
diffeomorphic to either
$S^4$ or $\bold R P^4$.  Assume in addition that 
\roster 
\item"{(i)}" $(M^4,g)$ is Bach-flat,
\smallskip
\item"{(ii)}" the Yamabe invariant $Y(M^4,g) > 0$,
\smallskip
\item"{(iii)}" the Weyl curvature satisfies
$$
\int_{M^4} |W|^2 dvol = 16 \pi^2 \chi(M^4).
\tag0.6
$$
\endroster
\smallskip\noindent
Then one of the following must be true:
\roster
\item"{1}" $(M^4,g)$ is conformal to $\bold C P^2$ with the 
Fubini-Study metric $g_{FS}$, or
\smallskip
\item"{2}" $(M^4,g)$ is conformal to a manifold which is 
isometrically covered by $S^3 \times S^1$
endowed with the product metric $g_{prod.}$.
\endroster
\endproclaim

In local coordinates, the Bach tensor is given by
$$
B_{ij} \, = \, 
\nabla^k 
\nabla^\ell 
W_{kij\ell} \, + \,
\frac{1}{2} \, 
R^{k \ell} \, W_{k i j \ell}.
\tag0.7
$$
Thus, Bach-flatness implies that the Weyl tensor lies in the kernel of a second order differential
operator.  At the same time, appealing once more to the results of \cite {CGY1} and \cite {CGY2},
we can prove that a manifold satisfying the hypotheses of Theorem C 
is conformal to one for which the integrand in (0.5) is identically zero. 
In section 4 we show how these facts,
along with a complicated Lagrange-multiplier argument, 
leads to the classification in Theorem C.

Actually, the proof of Theorem C is the most technically demanding 
aspect of the present paper.
First, there is a long calculation to derive an integral identity 
for the covariant derivative of the self-dual and anti-self-dual 
parts of the Weyl tensor of a Bach-flat metric.  
In addition to the algebraic difficulties of analyzing the curvature 
terms which arise in this identity, there are delicate
analytic issues.  For example, the conformal metric we construct 
based on the work of \cite {CGY1} and \cite {CGY2} may not be regular 
on all of $M^4$.  Indeed, if it were known to be smooth, 
then we could appeal to the classification of metrics with 
$WP \equiv {1 \over 6}$ done by Margerin.  These regularity 
problems are the price we pay, so to speak, for passing  
from integral to pointwise conditions.

We conclude the introduction with a note about the organization 
of the paper.  In Section 1 we develop the necessary PDE material 
from \cite {CGY1} and \cite {CGY2}.  Most of the
results are fairly straightforward generalizations of our earlier work. 
In Section 2 we show how the results of Section 1 and the work 
of Margerin can be combined to prove Theorem A.  Then, in Section 3 
we lay the groundwork for the proof of Theorems B and C by deriving 
various identities for the curvature of Bach-flat metrics.  In Section
4 we use these identities, along with an existence result from 
Section 1, to derive a key inequality for a certain 
polynomial in the curvature.  Analyzing this inequality leads 
us to consider a difficult Lagrange-multiplier problem, 
whose resolution gives the classification in the statement of Theorem C.

The research for this article was initiated while the
second author was a Visiting Professor at Princeton University
and the third author was a Visiting Member of the Institute for
Advanced Study, and was completed while all three authors
were visiting Institut des Hautes \'Etudes Scientifiques.  
The authors wish to acknowledge the support and 
hospitality of their host institutions.


\head{\bf 1.  An existence result}\endhead


In this section we prove an existence result in conformal geometry
which allows us to pass from the integral conditions of Theorems A--C 
to their pointwise counterparts.  As we indicated in the Introduction,
this result is based on the work in \cite{CGY1} and \cite{CGY2}.

To place this result in its proper context, we begin by introducing some
notation.  Given a Riemannian four--manifold $(M^4, g)$,
the {\it Weyl--Schouten} tensor is defined by
$$
A = Ric - {1\over 6} Rg
$$
In terms of the Weyl--Schouten tensor, the decomposition $(0.1)$
can be written
$$
\text{Riem} = W + {1\over 2} A \bcw g \,.
\tag1.0
$$
This splitting of the curvature tensor induces a splitting
of the Euler form. To describe this,
we introduce the elementary symmetric polynomials
$\sigma_\kappa : \Rn \rightarrow \R$,
$$
\sigma_\kappa ( \lambda_1 , \ldots , \lambda_n) =
   \sum\limits_{i_1 < \cdots < i_k}
     \lambda_{i_1} \cdots \lambda_{i_\kappa} \;.
$$
For a section $\Cal S$ of End$(TM^4)$ --- or, equivalently, a
section of $T^\star M^4 \otimes TM^4$ ---
the notation $\sigma_\kappa (\Cal S)$ means $\sigma_\kappa$ applied to the 
eigenvalues of $\Cal S$.  In particular, given a section of 
the bundle of symmetric two--tensors such as $A$, 
by ``raising an index'' we can cannonically associate a section
$g^{-1}A$ of  End$(TM^4)$.  At each point of $M^4$,
$g^{-1}A$ has 4 real eigenvalues, thus $\sigma_\kappa(g^{-1} A)$
is a smooth function on $M^4$.  To simplify notation, we
denote $\sigma_\kappa (A) = \sigma_\kappa (g^{-1}A)$.

Returning to the aforementioned splitting of the
Euler form, the Chern--Gauss--Bonnet formula (0.4) may be written
$$
8 \pi ^2 \chi (M^4) =
   \int {1\over 4 }|W|^2 dvol +
     \int \s2a dvol \;.
\tag1.1
$$
Note that the conformal invariance of the Weyl tensor implies
that the quantity
$$
\int \s2a dvol
$$
is conformally invariant as well.  Using (1.1),
assumption (0.5) of Theorem $A^\prime$
can be expressed
$$
\int \s2a dvol - \int {1\over 4} |W|^2 dvol > 0 \;.
\tag1.2
$$
Our goal in this section is to prove the following:

\proclaim{Theorem 1.1}
Let $(M^4 , g_0 )$ be a smooth, closed Riemannian
four--manifold for which 
\roster
\item"{(i)}"  The Yamabe invariant $Y(M^4 , g_0) >0$,
and
\item"{(ii)}" The curvature satisfies
$$
\int \sigma_2 (A_{g_0}) dvol - 
   \int {\alpha \over 4} | W_{g_0}|^2 dvol > 0 \,,
\tag1.3
$$
where $\alpha \geq 0$.  Then there is a conformal metric
$g_\alpha = e^{2w_{\alpha}} g_0$ whose curvature
satisfies
$$
\sigma_2 (A_{g_\alpha} ) - {\alpha\over 4} |W_{g_\alpha}|^2 \equiv \lambda \,,
\tag1.4
$$
where $\lambda$ is a positive constant.
\endroster
\endproclaim

Theorem 1.1 is a refinement of Corollary B of \cite{CGY2},
which for comparison's sake we now state:

\proclaim{Theorem 1.2}  (see \cite{CGY2}).
Let $(M^4, g_0)$ be a smooth, closed Riemannian
four--manifold for which
\roster
\item"{(i)}"  The Yamabe invariant $Y(M^4 , g_0) > 0$,
\item"{(ii)}"  The curvature satisfies
$$
\int \sigma_2 (A_{g_0} ) dvol_{g_0} > 0 \, .
\tag1.5
$$
\endroster
Then there is a conformal metric $g = e^{2w} g_0$
whose curvature satisfies
$$
\sigma_2 (A_g) \equiv \lambda \,,
\tag1.6
$$
where $\lambda$ is a positive constant.
\endproclaim

The proof of Theorem 1.2 relies on a crucial preliminary result:

\proclaim{Theorem 1.3} (See \cite{CGY1}).
Let $(M^4, g_0)$ be a smooth, closed Riemannian four--manifold
for which
\roster
\item"{(i)}"  The Yamabe invariant $Y(M^4, g_0) > 0$,
\item"{(ii)}"  The curvature satisfies
$$
\int \sigma_2 (A_{g_0}) dvol _{g_0} > 0.
\tag1.7
$$
\endroster
Then there is a conformal metric $g=e^{2w} g_0$ whose
curvature satisfies
$$
\sigma_2 (A_g) > 0 .
\tag1.8
$$
\endproclaim

The importance of Theorem 1.3 is that the metric satisfying inequality
(1.8) provides an approximate solution to equation (1.6), which
can then be deformed to an actual solution.  Moreover, (1.8)
implies that the path of equations connecting the metric
constructed in Theorem 1.3 to the metric constructed in
Theorem 1.2 is elliptic.

The proof of Theorem 1.1 will parallel the proofs of
Theorems 1.2 and 1.3.  The first step is the following
analogue of Theorem 1.3:

\proclaim{Theorem 1.4}
Let $(M^4 , g_0)$ be a smooth, closed Riemannian four--manifold
for which
\roster
\item"{(i)}"  The Yamabe invariant $Y(M^4 , g_0 ) > 0$,
\item"{(ii)}"  The curvature satisfies
$$
\int \sigma_2 (A_{g_0}) dvol_{g_0} -
  \int {\alpha\over 4} | W_{g_0}|^2
     dvol_{g_0} > 0 \, 
\tag1.9
$$
\endroster
where $\alpha \geq 0$.  Then there is a conformal metric
$g =e^{2w} g_0$ whose curvature satisfies
$$
\s2a - {\alpha\over 4} | W|^2 >0 .
\tag 1.10
$$
\endproclaim

\remark{Remark}
To simplify notation, we will use subscripts with $0$ 
instead of $g_0$, and
denote the volume form by $dv_0$ instead of $dvol_{g_0}$.
\endremark

\demo{Proof}
Following \cite{CGY1}, we consider the following functional
$F: W^{2,2} (M^4) \rightarrow \R$:
$$
F[\omega] =  \gamma_1 I[\omega] + \gamma_2 II [\omega] +
      \gamma_3 III [\omega]
$$
where $\gamma_i = \gamma_i(L)$ are constants and
$$
\align
I[\omega] &= \int 4|W_0|^2 \omega dv_0 -
    \left(\int |W_0|^2 dv_0 \right) \log \intll
       e^{4 \omega} dv_0 \,, \\
II[\omega] &= \int \omega P_0 \omega dv_0 +
    \int 4Q_0\omega dv_0 -
   \left( \int Q_0 dv_0 \right) \log \intll
       e^{4\omega} dv_0 \, , \\
III[\omega] &= 12 \left(Y[\omega] - {1\over 3}
   \int \Delta_0 R_0 \omega dv_0 \right) \,, \\
Y[\omega] &= \int (\Delta_0 \omega + |\nabla_0\omega|^2)^2
   dv_0 - {1\over 3 } \int R_0 | \nabla_0\omega|^2 dv_0\,.
\endalign
$$
Here $P$ denotes the Paneitz operator:
$$
P = (\Delta)^2 + d^\star \bigg({2\over 3} Rg - 2Ric \bigg)d\;,
$$
where $d$ is the exterior derivative , $d^\star$ is the adjoint of $d$,
and $Q$ is the fourth order curvature invariant:
$$
Q = {1\over 12} \left( - \Delta R + {1\over 4} R^2 - 3|E|^2 \right) \,.
$$
Thus
$$
Q = {1\over 2 }\s2a + {1\over 12} (- \Delta R) \,.
\tag1.11
$$
\enddemo

As in \cite{CGY2}, we need to introduce an additional
functional, which depends on the choice of a nowhere--vanishing
symmetric $(0,2)$--tensor $\eta$.
We then let
$$
\widetilde{I}[\omega] = \int 4|\eta|_0^2 \omega dv_0 -
    \left( \int | \eta|_0^2 dv_0 \right)
    \log \intll e^{4\omega} dv_0 \,.
$$
Now consider the functional
$$
F[\omega] = \widetilde{\gamma}_1 \widetilde{I} [\omega] +
 \gamma_1 I [\omega] + \gamma_2 II[\omega]
    +\gamma_3 III[\omega] \,,
$$
and define the conformal invariant
$$
\kappa = \widetilde{\gamma}_1 \int |\eta|_0^2 dv_0 +
   \gamma_1 \int | W_0|^2 dv_0 + \gamma_2
      \int Q_0 dv_0 \,.
\tag1.12
$$
Following the work of \cite{CY1}, we have
the following existence result for extremals of $F$. To make
the paper as self-contained
as possible, we will provide  
a sketch of the proof. 

\proclaim{Theorem 1.5}
(See \cite{CY1} Theorem 1.1)
Let $(M^4 , g_0)$ be a compact Riemannian four--manifold.
If $\gamma_2 , \gamma_3 >0$ and $\kappa < \gamma_2 8 \pi^2$,
then $\inf F[\omega]$ is attained by some function
$\omega \in W^{2,2} (M^4)$.  Moreover, the metric
$g = e^{2 \omega} g_0$ is smooth (see \cite{CGY3},\cite{UV})
and satisfies
$$
\widetilde{\gamma _1} | \eta |^2 + \gamma_1 |W|^2 +
   \gamma _2 Q - \gamma_3 \Delta R = \kappa vol (g)^{-1} \,.
\tag1.13
$$
\endproclaim

\demo{Proof}
To see that  $\inf F[\omega]$ is attained under the assumption that
$\gamma_2 , \gamma_3 >0$ and $\kappa < \gamma_2 8 \pi^2$, we employ a
sharp version of the Moser-Trudinger inequality established by
D.Adams [Ad]: there exists a constant $C= C(M, g_0)$ such that for
all $\omega \in W^{2.2}(M, g_0)$
$$
log \intll e^{ 4 (\omega  - \bar \omega )} dv_0 \leq C + \,\, \frac {1}{8 \pi^2}
\int (\Delta_0 \omega)^2 dv_0,
\tag 1.14
$$
where $ \bar \omega = \intll \omega dv_0 $.

Define
$$ 
U_0 = U(g_0) = \tilde \gamma_1 |\eta|_0^2 + \gamma_1 |W|_0 ^2 + \gamma Q_0
- \gamma_3 \Delta_0 R_0,
$$
then
$$
\int U_0 dv_0 = \kappa,
$$
and we can express $F$ as
$$ F[\omega] = - \kappa log \intll e^{ 4 (\omega - \bar \omega) } dv_0 + 4 \int U_0 ( \omega -
\bar \omega) dv_0 + \gamma_2 < P_0 \omega, \omega> + 12 \gamma_3 Y(\omega).
$$
When $\kappa \leq 0 $, we have
$$ F[\omega] \geq 4 \int U_0 ( \omega -
\bar \omega) dv_0 + \gamma_2 < P_0 \omega, \omega> + 12 \gamma_3 Y(\omega).
$$
When $ \kappa \geq 0$, we have from (1.14) that 
$$
\align
F[\omega] \geq & -C \kappa - \frac {\kappa}{ 8 \pi^2} \int (\Delta_0 \omega)^2 dv_0\\
+ &  4 \int U_0 ( \omega -
\bar \omega) dv_0 + \gamma_2 < P_0 \omega, \omega> + 12 \gamma_3 Y(\omega).
\endalign
$$
Thus, for a minimizing sequence $\{\omega_l\}$, $\underline {lim} F [\omega_l] \leq
F[0] = 0.$   From the estimates above we conclude that for $l$ large,
$\epsilon$ small,
$$
\align
\epsilon \geq F[\omega_l] \geq &  -C \tilde \kappa + ( - \frac {\tilde \kappa}{ 8 \pi^2} +
\gamma_2 + 12 \gamma_3 ) \int (\Delta_0 \omega_l)^2 \\ & + 12 \gamma_3 \int
|\nabla_0
\omega_l|^4 + 4 \int U_0 (\omega_l - \bar {\omega} _l )\\ & + ( \frac 23 \gamma_2 - 4
\gamma_3) \int R_0 |\nabla_0 \omega_l|^2 - 2 \gamma_2 \int Ric_0(\nabla_0 \omega_l,
\nabla_0 \omega_l)
\\ & + 24 \gamma_3 \int (\Delta_0 \omega_l) |\nabla_0 \omega_l|^2,
\endalign
$$
where $\tilde \kappa = max ( \kappa, 0)$.
It follows that if $ \kappa \leq \gamma_2 8 \pi^2$, and $\gamma_2 > 0$, $\gamma_3 >
0$, one concludes that there exists a
constant $C( g_0)$ so that
$$
\int (\Delta_0 \omega_l)^2 + |\nabla_0 \omega_l|^4 \leq C(g_0).
\tag 1.15
$$
Since the functional $F$ is scale-invariant, we may assume without
loss of generality that $ \intll \omega_l dv_0 = 0$. It follows from 
the Poincare inequality and (1.15) that $||\omega_l||_{ L^2} $ is uniformly bounded.
Therefore, 
$||\omega_l||_{2,2}$ is bounded and a subsequence will converge
weakly in $W^{2,2}$ to some $\omega \in W^{2,2}(M, g_0)$ with
$ F[\omega] = \inf _{\omega \in W^{2,2}} F[\omega]$.

Finally, a straightforward computation (see [BO]) shows that 
the extremal metric $g = e^{2\omega} g_0$
satisfies the Euler equation (1.13). 
\hskip 3.5in $\square$
\enddemo

Next, let $\delta \in (0,1]$ and choose
$$
\align
\widetilde{\gamma}_1 &= \left( -{1\over 2}\right)
  {\int \sigma_2 (A_0 ) dv_0 - {\alpha \over 4}
     \int | W_0 |^2 dv_0 \over \int |\eta|_0^2 dv_0} < 0 \,, \\
\gamma_1 &= - {\alpha\over 8} \, , \\
\gamma_2 &= 1 \,, \\
\gamma_3 &= {1\over 24 } (3\delta -2) \,.
\endalign
$$
With these values of $(\widetilde{\gamma}_1, \gamma_1, \gamma_2, \gamma_3$),
the conformal invariant $\kappa$ defined in (1.12) is equal to zero.
Thus, as long as $\gamma_3 > 0$ (i.e., $\delta > {2\over 3}$)
the hypotheses of Theorem 1.5 hold.  In particular, if
$\delta =1$ there exists an extremal metric
$g_1 = e^{2 \omega} g_0$, satisfying (1.13), which we rewrite using
(1.11):
$$%
\s2a - {\alpha\over 4} |W|^2 =
  {1\over 4} \Delta R - 2 \widetilde{\gamma}_1 | \eta|^2 \,.
\tag1.16
$$
Using the minimum principle of \cite{Gu, Lemma 1.2},
this implies that the scalar curvature of $g_1$
is strictly positive.

For general $\delta \in (0,1]$, the Euler equation (1.13) can be written
$$
\s2a -
{\alpha \over 4} |W|^2 = {\delta\over 4} \Delta R -
    2 \widetilde{\gamma}_1 | \eta|^2 \,.
\tag"{($\star\star)_\delta$}"
$$
Compare this with equation $(\star)_\delta$ of \cite{CGY1}:
$$
\s2a =
 {\delta\over 4} \Delta R - 2 \gamma_1 | \eta |^2 \,.
\tag"{($\star)_\delta$}"
$$
Note that the only discrepancy between $(\star\star)_\delta$
and $(\star)_\delta$ is the presence of the
Weyl term on the LHS of $(\star )_\delta$.
The key point is that this term has the same
sign as the term involving $\eta$
and scales in exactly the same way.  That is, if $g=e^{2w}g_0$ is
a conformal metric, then 
$$
\align
|W_g|_g^2 &= e^{-4w}|W_0|_0^2, \\
|\eta|_g^2 &= e^{-4w}|\eta|_0^2.
\endalign
$$
Consequently, the subsequent
arguments of section 2--6 of \cite{CGY1} can be carried
out with only trivial modifications, as we now describe.

Fixing $\delta_0 >0$, let
$$
\Cal S = \left \{ \delta \in [\delta_0 , 1] \quad\bigg| \;
   \foldedtext\foldedwidth{1.7in}{$(\star\star)_\delta$
     admits a solution with positive scalar cuvature}
\right\} \; .
$$
As we saw above, $ 1\in \Cal S$; thus $\Cal S$ is non--empty.
To verify that $\Cal S$ is open we compute the linearization of
$(\star\star)_\delta$, exactly as in Proposition 4.1 of \cite{CGY2}.
Again, the only relevant properties of the Weyl term and $\eta$--term
in $(\star\star)_\delta$ are their scaling properties 
(which are the same) and their sign (ditto).
Next, the estimates of section 3 in \cite{CGY1}
can be be used to show that $\Cal S$ is closed.
Consequently, for each $\delta >0$, there is a conformal
metric $g= g_\delta = e^{2\omega_\delta}g_0$ of
positive scalar curvature satisfying $(\star\star)_\delta$.

In sections 4--6 of \cite{CGY2} we obtain the following {\it a priori}
estimate for solutions of $(\star)_\delta$:  For fixed $p < 5$, there
is a constant $C= C(p)$ such that
$$
\Vert \omega_\delta \Vert _{2,p} \leq C \,.
\tag1.17
$$
In particular, $C$ is independent of $\delta$.
The same estimate holds for solutions of $(\star\star)_\delta$,
for the reasons explained above.

In section 7 of \cite{CGY1} the Yamabe flow is used to
show that one can perturb solutions of $(\star )_\delta$
to find metrics with $\s2a >0$.  An analogous result is true for
solutions of $(\star\star)_\delta$.

\proclaim{Theorem 1.6}
(See Theorem 7.1 of \cite{CGY1}).
Let $g = e^{2\omega} g_0$ 
be a solution of $(\star\star)_\delta$ with positive
scalar curvature, normalized so that $\int \omega dv_0 = 0$.
If $\delta >0$ is sufficiently small, then there is a smooth
conformal metric $h= e^{2v} g$ with
$$
\sigma_2(A_h) -{\alpha\over 4} |W_h|^2 > 0 \,.
\tag1.18
$$
\endproclaim

\demo{Proof}
The proof of Theorem 1.6, like the proof of its 
counterpart Theorem 7.1 of \cite{CGY1}, is based on
careful estimates of solutions to the Yamabe flow:
$$
\align
{\partial h\over \partial t} &= -{1\over 3} Rh , \\
h(0, \cdot) &= g = e^{2\omega} g_0 \,.
\endalign
$$
Using the estimate (1.17), we show that there
is a time $T_0$, which only depends on the background
metric $g_0$, such that the metric $h = h(T_0, \cdot)$
satisfies (1.18).  Although the arguments are essentially
the same, there are some necessary modifications of the
proof of Theorem 7.1 which require explanation.

First, Propositions 7.2 and 7.4, and Lemmas 7.3 and 7.5 can all
be copied without change.  In the statement of Proposition 7.8
we need to make an obvious change:  instead of defining
$f = \s2a + 2\gamma_1 | \eta |^2$, we define
$\tilde{f}=\s2a-{\alpha\over 4}|W|^2+2 \widetilde{\gamma}_1|\eta|^2$.
The conclusion of Proposition 7.8 then holds with
$f$ replaced by $\tilde{f}$.

In fact, by substituting $f$ with $\tilde{f}$,
the proof of the next proposition (7.12) is also valid.
Therefore, by following the remaining arguments of section 7
we arrive at the following inequality: For fixed
$s \in (4,5)$ and $t \leq T_1 (g_0)$,
$$
\s2a - {\alpha\over 4} |W|^2 \geq 
  -2 \widetilde{\gamma}_1 |\eta|^2 -C_3 t^{1-{4\over s}} 
   -C_3 \delta^{{1\over 2}} t^{(-1+{2\over s})}, 
$$
where $C_3 = C_3(g_0)$.  Since $s > 4$ and
$-2 \widetilde{\gamma}_1 |\eta | ^2 \geq C(g_0 ) > 0$,
it follows there is a constant
$C_4 = C_4 (g_0 ) > 0$ so that for
$t \leq T_0 = T_0 (g_0) $,
$$
\s2a - {\alpha \over 4} |W|^2 \geq
  {3\over 4} C_4 - C_3 \delta^{{1\over 2}} 
        t^{-(1+{2\over s})} \,.
\tag1.19
$$
Therefore, if $h = (t_0, \cdot )$, then for
$\delta < \delta _0 ( g_0 )$, (1.19) implies
$$
\sigma_2 (A_h) - {\alpha \over 4} |W_h|^2 \geq 
    {1\over 2} C_4 > 0 \,.
$$
This completes the proof of Theorem 1.6,
and consequently Theorem 1.4.  \hskip 1.09in $\square$

Once the existence of a metric satisfying (1.18)
is established, to complete the proof of
Theorem 1.1 we need to show that the techniques of \cite{CGY2} 
can be applied to construct a solution of (1.4).  This is actually
a two step process:  First, we need to establish {\it a priori}
bounds for solutions of
$$
\sigma_2 (A_g ) - {\alpha\over 4} |W_g|^2 = f > 0.
\tag1.20 
$$

The second step is to apply a degree-theoretic argument  
showing that a metric satisfying (1.18) can be deformed to a 
metric satisfying (1.4).  Of course, such an argument relies
on the estimates established in the first step.

\proclaim{Proposition 1.7}
(See Main Theorem of \cite{CGY2}).
Let $g = e^{2w} g_0$ 
be a solution of (1.20) with positive
scalar curvature, and assume $(M^4,g_0)$ is not conformally equivalent
to the round sphere.  
Then there is a constant $C = C(g_0,\|f\|_{C^2})$
such that
$$
\max_{M^4}\{e^w + |\nabla_0 w|\} \leq C.
\tag1.21
$$
\endproclaim

Now the estimate of \cite{CGY2}
applies to equations of the form
$$
\sigma_2(A_g) = f > 0,
$$
whereas (1.20) includes the Weyl term.  However, the argument of   
\cite{CGY2} can easily be modified to cover this case, as we now explain.

As in \cite{CGY2} we argue by contradicition: assuming the theorem
is false, then there is a sequence of solutions $\{w_{\kappa}\}$
of (1.20) (with $f$ fixed) such that
$$
\max_M [ | \trdo w_{\kappa} | + e^{w_{\kappa}} ] \rightarrow \infty
    \quad \text{ as }  \quad \kappa \rightarrow \infty \; .
$$
We then apply the blow-up argument described pages 155-156 of
\cite{CGY2}.  To begin, 
assume that $P_{\kappa} \in M^4$ is a point at which
$(|\nabla _0 w_{\kappa} | + e^{w_k})$ attains its
maximum.  By choosing normal coordinates $\{\Phi_{\kappa} \}$
centered at $P_{\kappa}$, we may identify the coordinate neighborhood
of $P_{\kappa}$ in $M^4$ with the unit ball $B(1) \subset \bold R^4$
such that $\Phi_{\kappa} (P_{\kappa}) = 0$.   Given $\varepsilon >0$,
we define the dilations $T_\varepsilon \colon \bold R^4 \rightarrow
\bold R^4$
by $x \mapsto \varepsilon x$, and consider the sequence
$w_{\kappa,\varepsilon} = T_\varepsilon^\star w_{\kappa} + \log \varepsilon$.
Note that
$$
| \trdo w_{\kappa,\varepsilon }| + e^{w_{\kappa,\varepsilon}}
  =\varepsilon (|\trdo w_{\kappa} | + e^{w_{\kappa}} ) \circ T_\varepsilon \; .
$$
Thus, for each $\kappa$ we can choose $\varepsilon_{\kappa}$ so that
$$
|\trdo (w_{\kappa, \varepsilon_{\kappa} }) | + e^{w_{\kappa,\varepsilon_{\kappa}}} 
      \big|_{x=0} = 1 \;.
$$
Note that $w_{\kappa, \varepsilon_{\kappa}}$ is defined in
$B_{\tsize{\frac{1}{ \varepsilon_{\kappa}}}} (0)$, and
$$
|\trdo (w_{\kappa, \varepsilon_{\kappa} )}| + e^{w_{\kappa}, \varepsilon_{\kappa}}
   \leqslant 1 \quad \text{ on } \quad
   B_{1\over \varepsilon_{\kappa}} (0) \; .
$$

To simplify notation, let us denote $w_{\kappa, \varepsilon_{\kappa}}$ by
$w_{\kappa}$.  Since from now on we view $\{w_{\kappa}\}$ as a sequence
defined on dilated balls in $\bold R^4$, there will be no
danger of confusing the renormalized sequence with the
original sequence.  Note that
$g_{\kappa}^\star \equiv e^{2w_{\kappa}} T_{\varepsilon_{\kappa}} ^\star g_0 \equiv e^{2w_{\kappa}} g_0^{\kappa}$
satisfies
$$
\sigma_2 (A_{g_{\kappa}^\star})
   - {\alpha \over 4} |W_{g_{\kappa}^\star}|^2 = f \circ T_{\varepsilon_k}  \; .
\tag1.22
$$
Furthermore, 
$g_0^{\kappa} = T_{\varepsilon_{\kappa}}^\star g_0$ $\rightarrow ds^2$;
where $ds^2$ is the Euclidean metric on $ \bold R^4$, in
$C^{2, \beta}$ on compact sets.  

As in \cite{CGY2}, we now have to consider two possibilities, depending on 
the behavior of the exponential term $e^{w_{\kappa}(0)}$.  However,
from here on the argument is identical in its details with that of
\cite{CGY2}.  The main point is that the conformal invariance of the
Weyl curvature implies that the two possible limit metrics 
arising from the sequence $g_{\kappa}^\star$ satisfy the same
equations as they do in \cite{CGY2}; that is, the Weyl term 
in (1.22) converges to
zero because Euclidean space is conformally flat.
In particular, Corollary 1.3 in 
\cite{CGY3} applies to equations such as (1.20), so the 
estimates needed to construct the limiting metric are the same.  

To summarize: after applying the same blow-up argument to our sequence, we
end up with the same limiting equations on
$\bold R^4$ (see Corollary 1.4 in \cite{CGY2}).  The rest of the
proof carries through exactly as in \cite{CGY2}, and we conclude that the
manifold $(M^4,g)$ must be conformally equivalent to the round 4-sphere.  Since
this contradicts our assumption, the estimate (1.21) must hold.

The local estimate of \cite{CGY2, Cor 1.3} also applies, and consequently we have a bound
$$
\Vert \nabla^2 w \Vert_\infty \leq C\,.
\tag 1.23
$$

Next, we use 
a degree theoretic argument
to prove the existence of a solution of (1.4). 
The following Proposition is 
a fairly straightforward modification of Corollary B in
\cite{CGY2}.

\proclaim {Proposition 1.8}
Assume that $(M^4, g_0)$ satisfies conditions (i) and (ii) of Theorem 1.4. 
Then given any positive (smooth) function $f \, > 0$, there exists a solution
$ g = e^{2w} g_0$ of (1.20). In particular, this is true if $ f \equiv \lambda$ for a constant $\lambda > 0$.
\endproclaim

\demo {Proof}
We first apply Theorem 1.4 to assert the existence of a conformal metric 
$g= e^{2w} g_0 $ for which the equation (1.20) holds for some positive function $f$. Since by assumption $(M^4, g_0)$ is not conformally equivalent
to the standard 4-sphere, the {\it a priori} estimates (1.22) and (1.23) 
hold.  In particular, given a smooth function $h$, there is a constant 
$c$ independent of $t$ so that all solutions
$g = e^{2w} g_0$ of the equation
$$
\sigma_2 (A_g ) - \frac {\alpha}{4} |W_g|^2 = t f + (1-t) h
\tag{$\Sigma_t$}
$$
with $R = R_g>0$ satisfy the bounds
$$
 \Vert w \Vert_{4, \alpha} \leqslant c, \quad
 S_{ij} (g) \xi_i \xi_j \geqslant \tfrac{1}{c} \, \vert \xi \vert
^2,
\tag 1.24
$$
where $S = -Ric + \frac{1}{2}Rg$ (see \cite{CGY1}, Lemma 1.2).
Let  $O_c$ be the set
$$
O_c = \{w \in C^{4,\alpha} : ( 1.24) \text{ holds } \}
   \cap \{ w \in C^{4,\alpha} | (\sigma_2(A_{g_w} - \frac {\alpha}{4} |W_{g_w}|^2 ) >0 ; \,\,
R_{g_w} > 0  \}.
$$
We denote the degree of the equation $(\Sigma_t)$ by $\deg
(\Sigma_t , O_c , 0)$.  The degree theory of [Li] implies that
$$
deg (\Sigma_0 , O_c , 0) = deg (\Sigma_1, O_c, 0 ).
\tag 1.25
$$
We need to do a calculation verifying that for $t=1$ the degree
of the equation is non-zero. In order to do this, we deform the
equation to one whose degree is easy to determine. First,
it is useful to re-write  equation (1.20) in
a suggestive form. Suppose $g= e^{2w} g_0$ and denote
$$
M_{ij}(w) = 2 S^{0}_{ij} + 2 \nabla^{0}_i\nabla^{0}_j w -
2\Delta_0 w g^{0}_{ij} -2 \nabla^{0}_i w \nabla^{0}_j w.
\tag 1.26
$$
Then, after some computation, the equation (1.20) may be written
in the form
$$
-\nabla^{0}_i\{M_{ij}(w)\nabla^{0}_j w\} + (\sigma _2 (A_{g_0}) - \frac{\alpha}{4} |W_{g_0}|^2) =
(\sigma_2 (A_g) - \frac {\alpha}{4} |W_g|^2) e^{4w} = \, fe^{4w}.
\tag 1.27
$$
It is important to note the identity
$$
M_{ij}(w) = S_{ij} + S^{0}_{ij} + |\nabla_{0} w|^2 g^{0}_{ij},
$$
so that it is clear that when both $ (\sigma_2 (A_g) - \frac {\alpha}{4} |W_g|^2) > 0 , R_g >0 $
and $ (\sigma_2(A_{g_0} - \frac {\alpha}{4} |W_{g_0} |^2)>0, R_{g_0}>0 $ , then $M_{ij}$ is
positive definite.

It is also convenient to re-formulate equation (1.20), on account
of the conformal covariance property, using the solution metric
$g$ of the equation as the background metric:
$$
-\nabla _i\{M_{ij}(v)\nabla _j v\} + f=f e^{4v},
\tag 1.28
$$
so that $v=0$ is a solution to this equation satisfying $R>0$.

We now use the following deformation:
$$
-\nabla _i\{M_{ij}(v)\nabla _j v\} + f= (\sigma_2(A_{g_v}) - \frac{\alpha}{4} |W_{g_v}|^2 ) e^{4v}
= (1-t)f\int e^{4v} + t f e^{4v}
$$
where $ \int e^{4v} = \int e^{4v} dvol_g $. We label this
equation by $\Gamma_t$. Note that when $t=1$, we recover the equation (1.28).
The proof of the Proposition now follows line by line the proof of
Corollary B in \cite{CGY2}.  More precisely, 
after establishing {\it a priori} estimates
for solutions of $\Gamma_t$, we find that the degree is well defined.
A calculation of the linearized equation,
together with the homotopy invariance of the degree implies that $deg(\Gamma_1, O_c, 0) = deg( \Gamma_0, O_c, 0) = -1$. It follows that a solution 
of (1.20) exists. 
\hskip 2.97in $\square$
\enddemo 

Applying Proposition 1.8, we complete the 
proof of Theorem 1.1. \hskip 2.97in $\square$
\enddemo


\head { \bf 2.   The proof of Theorem A}\endhead


Based on the results of section 1, we can now give
a detailed proof of Theorem A.

As we saw above, the assumption (0.3) is equivalent to the
inequality (1.2):
$$
\int\limits_{M^4} \s2a dvol -
   \int\limits_{M^4} {1\over 4} |W|^2 dvol > 0\,.
$$
Taking $\alpha =1$ in Theorem 1.1, it follows that
there is a conformal metric satisfying
$$
\s2a - {1\over 4} |W|^2 \equiv \lambda > 0\,.
\tag2.1
$$
Strictly speaking, at this stage all we really need is the conclusion
of Theorem 1.4 --- that is, we just need to know that the quantity 
in (2.1) is positive, not necessarily constant.

In any case, rewriting $\s2a$ in terms of the trace--free Ricci tensor
$E = Ric - {1\over 4} Rg$ and the scalar curvature we conclude
$$
- {1\over 2} |E|^2 + {1\over 24} R^2 -
    {1\over 4} |W|^2 > 0\,.
$$
Rearranging terms, this implies
$$
{|W|^2 + 2|E|^2 \over R^2} < {1\over 6} \,.
$$
By the weak-pinching result of Margerin \cite{Ma},
$M^4$ is diffeomorphic to $S^4$ or $\bold R P^4$. \hskip 1.1cm $\square$


\head{ \bf 3. A Weitzenb\"ock formula for Bach--flat metrics}\endhead


In preparation for the proof of Theorem C, in this section
we derive various curvature identities.  The first such result
is an inequality for metrics satisfying (1.4).

\proclaim{Lemma 3.1}
Suppose $(M^4, g)$ satisfies (1.4):
$$
\s2a - {\alpha \over 4} |W|^2 = \lambda\,,
$$
where $\alpha \geq 0$ and $\lambda \geq 0$ are constants.
Then
$$
{3\over 2 } \alpha |\nabla W|^2 + 3 ( |\nabla E|^2 -
   {1\over 12} | \nabla R|^2) \geq 0 \,.
\tag3.1
$$
\endproclaim

\demo{Proof}
Define the tensor
$$
V = \sqrt{\alpha} \,W + {1\over 2} E\bcw g\,,
$$
where $\bcw$ is the Kulkarni--Nomizu product.
Then
$$
\align
|V|^2 &= \alpha |W|^2 + 2 |E|^2
\tag3.2\\
|\nabla V|^2 & = \alpha |\nabla W|^2 + 2 | \nabla E|^2 \,.
\tag 3.3
\endalign
$$
As a consequence of (3.3), inequality (3.1) is
equivalent to
$$
|\nabla V|^2 \geq {1\over 6} |\nabla R|^2 \,.
\tag3.4
$$
To verify (3.4) note that (1.4) and (3.2) imply
$$
{1\over 6} R^2 = |V|^2 + 4 \lambda
\tag3.5
$$
Differentiating,
$$
{1\over 3} R \nabla R =
   \nabla |V|^2 = 2 |V| \nabla |V| \,.
$$
Taking the inner product of both sides with
$R^{-1} \nabla R$ gives
$$
{1\over 3} | \nabla R|^2 = 2 |V| 
   g\left(\nabla |V| , {\nabla R\over R}\right)
    \leq |V|^2 \, {|\nabla R|^2\over R^2} + 
       |\nabla |V||^2 \,.
$$
By Kato's inequality, $|\nabla |V||^2 \leq |\nabla V|^2$, and thus
$$
{1\over 3} | \nabla R|^2 \leq |V|^2 \,
{|\nabla R|^2\over R^2} + |\nabla V|^2 \,.
\tag 3.6
$$
Substituting (3.5) into (3.6) gives
$$
{1\over 3 }| \nabla R|^2 \leq \left(
 {1\over 6} R^2 - 4 \lambda \right) \,
   {|\nabla R|^2\over R^2} + |\nabla V|^2
  \leq {1\over 6} | \nabla R|^2 + | \nabla V|^2 \,,
$$
because $\lambda \geq 0$.  This establishes (3.4),
and consequently (3.1).\hskip 4.3cm $\square$
\enddemo

We now recall two important identities for Bach-flat metrics.
The first may be found in \cite{CGY1}:

\proclaim{Proposition 3.2}
(See \cite{CGY1}, Lemma 5.4).
If $(M^4 , g)$ is Bach-flat, then
$$
0 = \int\limits_{M^4} \left\{ 3 \, \left( |\nabla E|^2 - {1\over 12}
    |\nabla R|^2 \right) + 6tr E^3 + R|E|^2 -
      6W_{ijk\ell} E_{ik} E_{j \ell}\right\} dvol \,,
\tag3.7
$$
where $tr E^3 = E_{ij} E_{ik} E_{jk}$.
\endproclaim

\demo{Proof}
Just combine identities (5.5) and (5.10) of \cite{CGY1}. 
\qed
\enddemo

The second identity is a consequence of Stokes' Theorem,
the Bianchi identities, and the definition of the Bach tensor
in (0.7).

\proclaim{Proposition 3.3}
If $(M^4, g)$ is Bach flat, then
$$
\align
\hskip -1.6cm \int\limits_{M^4} |\nabla W |^2 dvol &=
  \int \left\{ 72 \det W^+ + 72 \det W^- -
  {1\over 2} R |W|^2 + 2 W _{ij\kappa\ell}E_{i\kappa} E_{j\ell} \right\}
\tag3.8
\endalign
$$
\endproclaim

\demo{Proof} 
In \cite{De}, Derdzinski proved a similar formula for metrics with 
harmonic Weyl tensor.  Our identity differs from his in only one
respect; namely, we replace harmonicity of the Weyl tensor 
(a first order condition) with Bach flatness (a second order condition).
To simplify the calculations and harmonize our notation with \cite{De},
we compute with respect to a local (normal) frame field.  Using this
convention, we write
$$
\int \bigl\vert \nabla W ^+ \bigr\vert^2 =
    \int \nabla_m W_{ijk\ell}^+ \nabla_m W_{ijk\ell} ^+ .
\tag3.9
$$
Since the splitting $\Lambda^2 = \Lambda_+^2 \oplus \Lambda_-^2$ is
parallel with respect to the Riemannian connection,
$$
\nabla W = \nabla W^+ + \nabla W^-. 
\tag3.10
$$
Therefore,
$$
\int \bigl\vert \nabla W^+ \bigr\vert^2 =
   \int \nabla_m W_{ijk\ell} \nabla_m W_{ijk\ell}^+ .
\tag3.11
$$
Using the decomposition (1.0), the second Bianchi
identity can be written
$$
\align
0 &= \nabla_m W_{ijk\ell} + \nabla_i W_{jmk\ell} + \nabla_j W_{mik\ell }\\
&\quad +{\tsize{1\over 2}} \bigl[ g_{ik} (dA)_{mj\ell}- g_{i\ell} (dA)_{mjk} \\
&\quad - g_{jk} (dA)_{mi\ell} + g_{j\ell} (dA)_{mik} - g_{km} (dA)_{ij\ell}\\
&\quad + g_{m \ell} (dA)_{ijk} \bigr] ,
\tag3.12
\endalign
$$
where
$$
(dA)_{ijk} = \nabla_i A_{jk} - \nabla_j A_{ik} .
$$
Contracting (3.12) we get the identity
$$
(\delta W )_{ij\ell} \equiv \nabla_m W_{ijm\ell} 
   = {\tsize{1\over 2}} (dA)_{ij\ell} ,
\tag3.13
$$
where $\delta$ denotes the divergence.  Subtituting
(3.12) into (3.11) we get
$$
\align
\int \bigl\vert \nabla W^+ \bigl\vert^2 &= \int \nabla_m 
     W_{ijk\ell}^+ \bigl\{ - \nabla_i W_{jmk\ell} \\
& \quad - \nabla_j W_{mik\ell} + {\tsize{1\over 2}} 
     \big[g_{km} (dA)_{ij\ell}
         - g_{m \ell} (dA)_{ijk}\big]\bigr\},
\endalign
$$
because all other terms vanish due to the symmetries of the Weyl tensor.
Re--indexing and combining like terms, we find
$$
\int \big\vert \nabla W^+\big|^2 =
    \int -2 \nabla_m W_{ijk\ell}^+ \nabla_i W_{jmk\ell}
         +\nabla_m W_{ijm\ell}^+ (dA)_{ij\ell} .
\tag3.14
$$
The splitting (3.10) implies $\delta W = \delta W^+ + \delta W^- $,
so by (3.13)
$$
\nabla_m W_{ijm\ell} (dA)_{ij\ell} = 
     2(\delta W^+)_{ij\ell} (\delta W)_{ij\ell} = 2 |\delta W^+|^2 .
$$
Similarly,
$$
-2 \nabla_m W_{ijk\ell}^+ \nabla_i W_{jmk\ell} =
     -2 \nabla_m W_{ijk\ell}^+ \nabla_i W_{jmk\ell}^+ .
$$
Substituting these into (3.14) we obtain
$$
\int \big\vert \nabla W^+\big\vert ^2 =
   \int -2 \nabla_m W_{ijk\ell}^+ \nabla_i W_{jmk\ell}^+
      +2 \big\vert \delta W^+\big\vert^2 .
\tag3.15
$$
We analyze each term in (3.15) separately.

For the first term, we integrate by parts and commute derivatives:
$$
\align
\int -2 \nabla_m W_{ijk\ell}^+ \nabla_i W_{jmk\ell}^+ 
&= \int 2 W_{ijk\ell}^+ \nabla_m \nabla_i W_{jmk\ell}^+  \\
& = \int 2 W_{ijk\ell}^+ \bigl\{ \nabla_i \nabla_m W_{jmk\ell}^+ 
         +R_{mijs} W_{smk\ell}^+\\
&\qquad  + R_{mims}W_{jsk\ell}^+  + R_{miks} W_{jms\ell}^+ \\ 
& \qquad +   R_{mi\ell s}W_{jmks}^+ \bigr\} \\
& = \int 2 W_{ijk\ell}^+ \nabla_i \nabla_m W_{jmk\ell}^+
     + 2R_{mijs} W_{ijk\ell}^+ W_{smk\ell}^+ \\
&\qquad + 2R_{is}W_{ijk\ell}^+ W_{jsk\ell}^+
   +2R_{miks}W_{ijk\ell}^+ W_{jms\ell}^+  \\
& \qquad   +2R_{mi\ell s} W_{ijk\ell}^+ W_{jmks}^+ ,
\tag3.16
\endalign
$$
where $R_{is} = R_{mims}$ are the components of the Ricci tensor.
Note by re--indexing
the last two terms in (3.16) are equal.  If we integrate by parts again,
the first term in (3.16) is
$$
\int 2W_{ijk\ell}^+ \nabla_i \nabla_m W_{jmk\ell}^+
    = \int -2 \nabla_i W_{ijk\ell}^+ \nabla_m W_{jmk\ell}^+
      = \int 2\big\vert \delta W^+\big\vert^2 .
\tag3.17
$$
Using the Bianchi identity and re--indexing
the next term can be rewritten
$$
2R_{mijs} W_{ijk\ell}^+ W_{smk\ell}^+ 
    = R_{msij} W_{ijk\ell}^+ W_{msk\ell}^+ \,.
\tag3.18
$$
Appealing to the decomposition (1.0) once more, we get
$$
2R_{mijs} W_{ijk\ell}^+ W_{smk\ell}^+
   = W_{msij}^+ W_{ijk\ell}^+ W_{msk\ell}^+
     + 2A_{js} W_{jik\ell}^+ W_{sik\ell}^+ \, .
\tag3.19
$$
Similarly,
$$
4R_{miks} W_{ijk\ell}^+ W_{jms\ell}^+ 
    = 4W_{miks} W_{ijk\ell}^+ W_{jms\ell}^+
      +2A_{is} W_{ijk\ell}^+ W_{jks\ell}^+
      + 2A_{km} W_{ijk\ell}^+ W_{jmi\ell}^+ \, .
\tag3.20
$$
Using the symmetries of the Weyl tensor and re--indexing we find
$$
2A_{is} W_{ijk\ell}^+ W_{jks\ell}^+ + 2A_{km} W_{ijk\ell}^+ W_{jmil}^+
    = -2A_{km} W_{ij\ell k}^+ W_{ij\ell m}^+ \, .
\tag3.21
$$
Combining (3.15)--(3.21),
$$
\align
\int \big\vert \nabla W^+\big\vert^2 &=
    \int 4 \big\vert \delta W^+ \big\vert^2
      + 2R_{is} W_{ijk\ell}^+ W_{jsk\ell}^+ \\
&\qquad + W_{msij}^+ W_{ijk\ell}^+ W_{msk\ell}^+ 
        + 4 W_{miks} W_{ijk\ell}^+ W_{jms\ell}^+ \, .
\tag3.22
\endalign
$$

\proclaim{Lemma 3.4}
\roster
\item"{(i)}"  $W_{ijk\ell}^+ W_{jsk\ell}^+$
$ = - \, {1\over 4} \big\vert W^+\big\vert^2 \delta_{ij}$,
\item"{(ii)}"  $W_{msij}^+ W_{ijk\ell}^+ W_{msk\ell}^+ = 24 \det W^+ $,
\item"{(iii)}"  $4W_{miks} W_{ijk\ell}^+ W_{jms\ell}^+ = 48 \det W^+$.
\endroster
\endproclaim

\demo{Proof}
As in \cite{De}, we fix a point and diagonalize
$W^{\pm} : \Lambda_{\pm}^2 \longrightarrow \Lambda_{\pm}^2$.
Let $\lambda_i ^{\pm}, 1\leq i\leq 3$ denote the three
eigenvalues of $W^\pm$, with corresponding eigenforms
$\omega^\pm , \eta^\pm , \theta^\pm$.  Then
$$
W^\pm = \lambda_1^\pm \omega^\pm \otimes \omega^\pm +
   \lambda_2^\pm \eta^\pm \otimes \eta^\pm + 
       \lambda_3^ \pm \theta^\pm \otimes \theta^\pm .
$$
Then (i)--(iii) follow from elementry calculations.
\qed
\enddemo

From the preceding Lemma and (3.22) we obtain the identity
$$
\int \big\vert \nabla W^+\big\vert^2 =
    \int 4 \big\vert \delta W^+ \big\vert^2 
      +72 \det W^+ - {\tsize{1\over 2}} 
         R\big\vert W^+ \big\vert^2 ,
\tag3.23
$$
which holds for any Riemannian four--manifold.

Now suppose $(M^4, g)$ is Bach--flat.  By \cite{De, (23)},
$$
0= \nabla_k \nabla_\ell W_{ik\ell j} - {\tsize{1\over 2}}
        W_{ikj\ell} A_{k\ell} \, .
$$
Pairing both sides with the Weyl--Schouten tensor and
integrating we get
$$
\align
0&= \int A_{ij} \nabla_k \nabla_\ell W_{ik\ell j} 
     - {\tsize {1\over 2}} W_{ikj\ell}  A_{k\ell} A_{ij}\\
&= \int - \nabla_k A_{ij} \nabla_\ell W_{ik\ell j} 
      -{\tsize{1\over 2}} W_{ikj\ell} A_{k\ell} A_{ij} \\
&= \int - {\tsize{1\over 2}} \big(\nabla_k A_{ij} 
        - \nabla_i A_{kj} \big)
         \nabla_\ell W_{ik\ell j} - {\tsize{1\over 2}}
           W_{ikj\ell} A_{k\ell} A_{ij} \\
&= \int - {\tsize{1\over 2}} (dA)_{kij} \nabla_\ell W_{ik\ell j}
        - {\tsize{1\over 2}} W_{ikj\ell} A_{k\ell} A_{ij} \\
&= \int \big\vert \delta W\big\vert^2 
         - {\tsize{1\over 2}} W_{ikj\ell} A_{k\ell} A_{ij} \,.
\endalign
$$
Since the Weyl tensor is trace--free, we conclude
$$
\int \big\vert \delta W \big\vert^2 =
    \int {\tsize{1\over 2 }} W_{ijk\ell} E_{ik} E_{j\ell} .
\tag3.24
$$
Combining (3.23) and (3.24), we get
$$
\align
\int \big\vert \nabla W\big\vert^2 
 &  =  \int \big\vert\nabla W^+\big\vert^2 + 
         \int \big\vert \nabla W^-\big\vert^2 \\
&= \int 72 \det W^+ + 72 \det W^-
    - {\tsize{1\over 2}} R \big\vert W \big\vert ^2 +
      2 W_{ijk\ell} E_{ik} E_{j\ell} .
\endalign
$$
This completes the proof of Proposition 3.3.
\qed
\enddemo

Combining (3.1), (3.7) and (3.8) we find that for any
$\alpha \geq 0$,
$$
\align
0& = \int\limits_{M^4} \left\{ {3\over 2} \alpha
   |\nabla W|^2 + 3 \left(|\nabla E|^2 - {1\over 12} 
       |\nabla R|^2 \right) \right.\\
&\qquad\qquad + 6trE^3 + R|E|^2 - 3(\alpha +2) 
      W_{ij\kappa \ell} E_{i\kappa}E_{j\ell}  \\
&\qquad \qquad \left. - 108 \alpha \det W^+ - 108 \alpha \det W^- +
   {3\over 4} \alpha R | W |^2 \right\} dvol \\
&\geq \int\limits_{M^4} 6trE^3 + R|E|^2 - 3(\alpha +2) 
      W_{ij\kappa \ell} E_{i\kappa}E_{j\ell}  \\
&\qquad \qquad  - 108 \alpha \det W^+ - 108 \alpha \det W^- +
   {3\over 4} \alpha R | W |^2  \,.
\tag3.25
\endalign
$$

This is the key identity in the proof of Theorem C.


\head{ \bf 4.  The proof of Theorem C}\endhead


Suppose $(M^4, g)$ is Bach--flat with positive Yamabe
invariant, and that (0.6) holds:
$$
\int\limits_{M^4} |W|^2 dvol = 16 \pi^2 \chi (M^4) \,.
$$
This is equivalent to
$$
\int\limits_{M^4} \s2a dvol = {1\over 4} \int\limits_{M^4}
   |W|^2 dvol \,.
\tag4.1
$$
Now, if $W \equiv 0$ then $(M^4, g)$ is locally conformally
flat and (by (0.6)) $\chi(M^4) =0$.
It follows from \cite{Gu, Corollary G} that
$(M^4, g)$ is conformal to a manifold which is
isometrically covered by $S^3 \times S^1$.  Therefore,
let us assume from now on that
$$
\int\limits_{M^4} | W|^2 dvol > 0 \,.
\tag4.2
$$
By (4.1), this implies that for any $0\leq \alpha < 1$,
$$
\int\limits_{M^4} \s2a dvol - {\alpha \over 4}
   \int\limits_{M^4} |W|^2 dvol =
   {(1-\alpha )\over 4}
    \int\limits_{M^4} |W|^2 dvol > 0\,.
$$
According to Theorem 1.1, there is a conformal metric 
$g_\alpha = e^{2\omega_\alpha} g$ satisfying
$$
\sigma_2 (A_{g_\alpha}) - {\alpha\over 4}
   |W_{g_\alpha}|^2 \equiv \lambda_\alpha > 0 \,.
\tag4.3
$$
Choose a sequence $\alpha_\kappa \nearrow 1$, and 
denote $g_\kappa = g_{\alpha_\kappa} = e^{2\omega_\kappa} g$.
The compactness properties of the sequence $\{g_\kappa\}$ require
careful description.

First, we claim that the a priori estimate of
\cite{CGY2} holds: i.e., there is a constant $C$ such that
$$
\Vert \nabla \omega_\kappa \Vert_\infty +
  \Vert \omega_\kappa \Vert _ \infty \leq C \,.
\tag4.4
$$
(Recall (4.2) implies that $(M^4 , g)$ cannot be
conformally equivalent to the sphere).  This estimate is
immediate from Proposition 1.7.  In addition, 
the local estimate of \cite{CGY2, Cor 1.3} applies, and consequently we have a bound
$$
\Vert \nabla^2 \omega_\kappa \Vert_\infty \leq C\,.
\tag4.5
$$
It is important to note that (4.5) is optimal:
(4.3) is elliptic if $\s2a > 0$; but $\lambda_\kappa \rightarrow 0$
as $\kappa \rightarrow \infty$, and there is no guarantee that the
Weyl tensor does not vanish on $M^4$.  However, higher order
estimates for $\{g_\kappa\}$ can be established on the set
where the Weyl tensor is non--zero.
To explain this, let
$$
\align
M_+^4 & = \{ x \in M^4 : |W_g| > 0 \} , \\
M_0^4 &= \{ x \in M^4 : | W_g| = 0 \} \,.
\endalign
$$
By conformal invariance, the
Weyl tensor of each $g_\kappa$ is also non-zero on $M_+^4$ (and vanishes
on $M_0^4$).  If $x_0 \in M_+^4$, then there are constants
$\epsilon > 0 , \rho > 0$, such that $|W_g| \geq \epsilon >0$
on the geodesic ball 
$B_\rho (x_0 )=\{x \in M^4: dist_g(x,x_0)<\rho \}$.
On $B_\rho (x_0)$ the metric
$g_\kappa = e^{2 \omega \kappa} g$ satisfies
$$
\sigma_2(A_{g_\kappa}) = {\alpha_\kappa\over 4}
   |W_{g_\kappa} |^2 + \lambda_\kappa
   = {\alpha_\kappa\over 4} e^{-4 \omega_\kappa} | W_g|^2 + \lambda_\kappa \,.
$$
In particular, by the a priori estimate (4.4) we see that
$$
\sigma_2 (A_{g_\kappa}) \geq C_\epsilon > 0
$$
on $B_\rho (x_0)$.  Therefore,
$$
\sqrt{\sigma_2 (A_{g_\alpha})} = \sqrt{{\alpha\over 4}
   |W_{g_\alpha}|^2 + \lambda_\alpha} 
$$
is a strictly elliptic,
concave equation on $B_\rho (x_0)$.
The regularity results of Evans \cite{Ev}
and Krylov \cite{Kr} then give H\"older estimates for
$\nabla^2 \omega_\kappa$ on $B_{\rho^\prime}(x_0)$ for any $\rho^\prime < \rho$
(see \cite{GT, Theorem 17.14}).  Applying Schauder theory and
classic elliptic regularity, we then obtain estimates for
derivatives of all orders on any ball
$B_{\rho^\prime} (x_0) \varsubsetneq B_\rho (x_0)$
with $\rho^\prime < \rho$.
Consequently, a subsequence of $\{g_\kappa\}$ (also denoted  $\{g_\kappa\}$)
 converges to a limiting metric
$g_\infty$ in $C_{loc}^\infty (M_+^4) \cap C^{1,1}(M^4)$.

Recall that inequality
(3.25) is satisfied by each $g_\kappa$ (with $\alpha = \alpha_\kappa$).  If we split
the integral in (3.25) into two integrals, one over 
 $M_+^4$, the other over $M_0^4$, note that on
$M_0^4$ the integrand reduces to
$$
6 trE^3 + R|E|^2 \,.
$$
Using the sharp inequality
$trE^3 \geq - {1\over \sqrt{3}}\, |E|^3$, we find that
$$
6 trE^3 + R|E|^2 \geq - 2 \sqrt{3} \,|E|^2 + R|E|^2
 = |E|^2 \big(R - 2 \sqrt{3} \, |E| \big) \;.
\tag4.6
$$
Since $\s2a \geq 0$ on $M_0^4$, we also have
$$
\alignat2
&\qquad  \;0 \leq  - {\tsize{1\over 2}} |E|^2 + {\tsize{1\over 24}} R^2 \\
&\Rightarrow  R^2 \geq  12 |E|^2 \\
&\Rightarrow  \; R \geq 2 \sqrt{3}\,|E| \,.
\tag4.7
\endalignat
$$
Combining (4.6) and (4.7) we see that the integrand in (3.25) is non-negative
on the set $M_0^4$.  Therefore, in view of the convergence of 
$\{g_\kappa\}$ on $M_+^4$ we have 
$$
\align
0&\geq \lim_{\kappa \to \infty } \int\limits_{M_{+}^4} 6trE^3 + R|E|^2 - 3(\alpha_{\kappa} +2) 
      W_{ij\kappa \ell} E_{i\kappa}E_{j\ell}  \\
&\qquad \qquad  - 108 \alpha_{\kappa} \det W^+ - 108 \alpha_{\kappa} \det W^- +
   {3\over 4} \alpha_{\kappa} R | W |^2 \\
& =  \int\limits_{M_{+}^4} 6trE^3 + R|E|^2 - 9 
      W_{ij\kappa \ell} E_{i\kappa}E_{j\ell}  \\
&\qquad \qquad  - 108 \det W^+ - 108 \det W^- +
   {3\over 4} R | W |^2 \,.
\tag4.8
\endalign
$$

To summarize:  we have constructed a metric
$g_\infty = e^{2w_{\infty}}g$ on $M^4$ with $w_{\infty} \in C^\infty (M_+^4) \cap C^{1,1}(M^4)$
which satisfies
\CenteredTagsOnSplits
$$
\left.
\matrix
\hskip -1.4in \sigma_2(A) = \tsize{{1\over 4 }}|W|^2 \, \text{on } \, M_+^4, \\
\quad 0\geq \dsize{\int\limits_{M_+^4}} \Bigl\{ 
      -108 \det W^+ - 108 \det W^- \bigg. \\ 
\qquad \quad \qquad \qquad + 6trE^3 -  
       9 W_{ij\kappa \ell} E_{i\kappa} E_{j\ell}\\
\\
 \qquad \qquad \qquad \quad  \; \; \Bigl. 
       +R\Bigl(|E|^2 + \tsize{{3\over 4}} |W|^2 \Bigr) \Bigr\} dvol \,.
\endmatrix
\quad \right\}
\tag4.9
$$

\proclaim{Proposition 4.1}
If $g_\infty$ satisfies (4.9), then
\roster
\item"{(i)}" $g_\infty \in C^\infty (M^4)$;
\item"{(ii)}" $g_\infty$ is Einstein;
\item"{(iii)}" either $W^+ \equiv 0$ or $W^- \equiv 0$ on $M^4$.
\endroster
\endproclaim

\demo{Proof}
Our first task is to rewrite the integrand in (4.9) in a
suitable basis.
To this end, let $Riem: \Lambda^2 \rightarrow \Lambda^2$
denote the curvature operator of $(M^4 , g_\infty )$.
Since $M^4$ is oriented, we have the splitting 
$\Lambda ^2 = \Lambda_+^2 \oplus \Lambda_-^2$, and the well
known decomposition of Singer--Thorpe \cite{ST}:
$$
Riem = \left(
\matrix
W^+ + {1\over 6} R \, Id &\vdots& B \\
\dotfill &\vdots&\dotfill \\
B^\star &\vdots &W^- + {1\over 6} R \, Id
\endmatrix
\right)
\tag4.10
$$
\enddemo
Note the compositions satisfy
$$
\align
&BB^\star : \Lambda_+^2 \rightarrow \Lambda_+^2 \,, \\
&B^\star B : \Lambda_-^2 \rightarrow \Lambda_-^2 \, .
\tag4.11
\endalign
$$

Fix a point $P \in M_+^4$, and let
$\lambda_1^\pm \leq \lambda_2^\pm \leq \lambda_3^\pm$ denote
the eigenvalues of $W^\pm$.  Then
$$
\align
\det W^\pm & = \la123
\tag4.12 \\
|W^\pm |^2 &= 4\Vert W^\pm \Vert ^2 = 4[(\lambda_1^\pm)^2 +
    (\lambda_2^\pm )^2 + (\lambda_3^\pm )^2] \,.
\tag4.13 
\endalign
$$
Recall that $\Vert W^\pm \Vert$ denotes the norm of $W^\pm$
when interpreted as an endomorphism of $\Lambda_\pm^2$.

Following Margerin, we denote the eigenvalues of
$BB^\star : \Lambda_+^2 \rightarrow \Lambda_+^2$
by $\b12b22b32$, where $0 \leq \bB1bB2bB3$.

\proclaim{Lemma 4.2}
$$
\align
|E|^2 &= 4(b_1^2 + b_2^2 + b_3^2),
\tag4.14 \\
trE^3 &\geq - 24 b_1b_2b_3 \,.
\tag4.15
\endalign
$$
\endproclaim

\demo{Proof}
Given a basis $\{ e_i\}$ of $T_pM^4$,
let $\{e_i^\star \}$ denote the dual basis of $T_p^\star M^4$.
Relative to this basis, the curvature operator is given by
$$
Riem (e_i^\star \wedge e_j^\star ) =
{1\over 2} \sum\limits_{\kappa , \ell}
  R_{ij\kappa \ell} e_\kappa ^\star \wedge e_\ell ^\star \,,
\tag4.16
$$
where $R_{ij\kappa \ell}$ are components of $Riem$ viewed as a
$(0,4)$--tensor;  i.e., $R_{ij\kappa \ell} = Riem (e_i, e_j, e_\kappa, e_\ell)$.
If $E_{ij} = E(e_i, e_j )$ are the components of the trace--free
Ricci tensor, then the decomposition $(0.1)$ implies
$$
\split
R_{ij\kappa \ell} = W_{ij\kappa \ell} +
   {\tsize{1\over 2}} (\delta_{i \kappa} E_{j\ell} - \delta_{i\ell} E_{j \kappa}
       - \delta_{j\kappa} E_{i \ell}
  + \delta_{j\ell} E_{i\kappa} )\\
  + {\tsize{1\over 12}} R(\delta_{i\kappa} \delta_{j \ell} - \delta_{i \ell}
           \delta_{j \kappa }) \,,
\endsplit
\tag4.17
$$
where $W_{ij\kappa \ell} = W(e_i , e_j , e_\kappa , e_\ell )$.

A basis $\{e_i^\star \}$ of $T_p^\star M^4$ induces a natural
orthonormal basis of $\Lambda_\pm^2$:
$$
\align
\omega^\pm & =  {\tsize{1\over \sqrt{2}}} 
    (e_1^\star \wedge e_2^\star \pm e_3^\star \wedge e_4^\star )\,,\\
\eta^\pm &= {\tsize{1\over \sqrt{2}}} (e_1^\star 
      \wedge e_3^\star \mp e_2^\star \wedge e_4^\star )\,,\\
\theta^\pm &= {\tsize{1\over \sqrt{2}}} 
      (e_1^\star \wedge e_4^\star \pm e_2^\star \wedge e_3^\star )\,.\\
\tag4.18
\endalign
$$
Now suppose the basis $\{ e_i\}$ diagonalizes $E$:
$$
E= \pmatrix
E_{11}&&&&0\\
&E_{22}&&& \\
&&E_{33}&&\\
0&&&E_{44}& \\
\endpmatrix.
$$
Using (4.11), (4.16), and (4.17), the matrix of
$B: \Lambda_{-}^2 \rightarrow \Lambda_+^2$ relative to the basis in (4.18) is
$$
B=
\pmatrix
{\tsize{1\over 4}} (E_{11}+E_{22} -E_{33} - E_{44})&&0 \\
\\
&{\tsize{1\over 4}} (E_{11} + E_{33} - E_{22} -E_{44})& \\
\\
0&&{\tsize{1\over 4}} (E_{11} + E_{44} - E_{22} - E_{33}) \\
\endpmatrix.
$$
Let
$$
\align
\mu_1 &= {\tsize{1\over 4}} (E_{11} + E_{22} - E_{33} - E_{44}),\\
\mu_2 &= {\tsize{1\over 4}} (E_{11} + E_{33} - E_{22} - E_{44}),\\
\mu_3 & ={\tsize{1\over 4}} (E_{11}+ E_{44} - E_{22} - E_{33}) \,.
\endalign 
$$
Since $E$ is trace--free, these can also be expressed
$$
\align
\mu_1 &= {\tsize{1\over 2}} (E_{11} + E_{22}) , \\
\mu_2 &= {\tsize{1\over 2}} (E_{11} + E_{33}) , \\
\mu_3 &= {\tsize{1\over 2}} (E_{11} + E_{44}) .
\tag4.19
\endalign
$$
Consequently, in terms of $\{\mu_1 , \mu_2, \mu_3 \}$
the eigenvalues of $BB^\star: \Lambda_+^2 \rightarrow \Lambda_+^2$ are
$$
\align
b_1^2 = \mu_1^2 , \\
b_2^2 = \mu_2^2 , \\
b_3^2 = \mu_3^2 , \\
\tag4.20
\endalign
$$
Now, a simple calculation gives
$$
\align
8\mu_1\mu_2\mu_3&= (E_{11} + E_{22})(E_{11} + E_{33}) (E_{11}+E_{44})\\
&= E_{11}^3 + E_{11}^2E_{22} + E_{11}^2 E_{33}+E_{11}^2 E_{44}\\
&\qquad + E_{11}E_{22} E_{33} + E_{11}E_{22}E_{44} +E_{11} E_{33}E_{44} 
      + E_{22}E_{33}E_{44}\\
&= E_{11}^2 (E_{11} +E_{22} + E_{33} +E_{44}) +
     \sigma_3 (E_{11} , E_{22}, E_{33}, E_{44} ) .
\endalign
$$
On the other hand, for a symmetric trace--free
$4\times 4$ matrix $E$, $\sigma_3 (E) = {1\over 3} tr E^3$.  Thus
$$
\align
tr E^3 &= 24 \mu_1\mu_2\mu_3 \\
& \geq -24 |\mu_1 \mu_2 \mu_3 | \\
&= -24 b_1b_2b_3 \,.
\endalign
$$
This proves (4.15).  The proof of (4.14) follows from (4.20), 
and will be omitted. \hskip 1.3in $\square$
\enddemo

The next inequality follows from Lemma 6 in \cite{Ma2}.
However, as our notation and conventions are slightly different
we provide some details.  

\proclaim{Lemma 4.3}
$$
- W_{ij\kappa \ell} E_{i \kappa} E_{j \ell} \geq
-4 \left[ \sum\limits_{i=1}^3 \lambda_i^+ b_i^2 +
     \sum\limits_{i=1}^3 \lambda_i^- b_i^2 \right]
\tag4.21
$$
\endproclaim

\demo{Proof}
This inequality is termed ``decoupling of the Weyl and Ricci curvatures''
by Margerin, and appropriately enough:  In general $W^\pm$
and $BB^\star$ or $B^\star B$ do not commute, and therefore cannot
be simultaneously diagonalized.

In any case, if we choose a basis $\{e_i\}$ of $T_pM^4$ which diagonalizes
$E$ as in Lemma 1, then the matrix of 
$W^\pm : \Lambda_\pm ^2 \rightarrow \Lambda_\pm^2$
relative to the basis in (4.18) is
$$
W^\pm =
\pmatrix
{1\over 2}(W_{1212}\pm 2W_{1234} +W_{3434})&&\star\\
&{1\over 2}(W_{1313}\mp 2W_{1324} +W_{2424})&\\
\star &&{1\over 2}(W_{1414}\pm W_{1423} +W_{2323}) \\
\endpmatrix .
$$
Therefore, if $\langle \quad , \quad \rangle_{\Lambda_{\pm}^2}$
denotes the natural inner product induced on $\Lambda_\pm ^2$
(see Remark 2 in the Introduction), then
$$
\align
\langle W^+ , BB^\star \rangle_{\Lambda_{+}^2} + 
\langle W^- , B^\star B \rangle_{\Lambda_{-}^2} 
  &= tr (W^+ \circ BB^\star + W^- \circ B^\star B) \\
&  = \left\{  (W_{1212} + W_{3434}) \mu_1^2 + 
              (W_{1313} + W_{2424})\mu_2^2 \right.\\
& \qquad +  \left.   (W_{1414} + W_{2323} ) \mu_3^2 \right\}\\
&= {\tsize{1\over 2}}(W_{1212} E_{11} E_{22} + W_{1313} E_{11} E_{33}
     +W_{1414}E_{11}E_{44} \\
&\qquad +W_{2323} E_{22}E_{33} +W_{2424}E_{22}E_{44}
         +W_{3434}E_{33} E_{44} ) .
\endalign
$$
On the other hand,
$$
\align
W_{ij\kappa \ell}E_{i\kappa}E_{j\ell} &=
   2(W_{1212}E_{11}E_{22} + W_{1313} E_{11}E_{33}
     +W_{1414}E_{11}E_{44} \\
  & \qquad + W_{2323}E_{22} E_{33}
     +W_{2424}E_{22}E_{44} + W_{3434} E_{33} E_{44} ).
\endalign
$$
Therefore,
$$
-W_{ij\kappa \ell} E_{i\kappa} E_{j \ell} =
   -4\langle W^+ , BB^\star \rangle_{\Lambda_+^2}
     -4 \langle W_1^- B^\star B\rangle_{\Lambda_-^2} \;.
\tag4.22
$$
According to Lemma 6 of \cite{Ma},
$$
\langle W^+, BB^\star \rangle_{\Lambda_+^2}  +
    \langle W^- , B^\star B \rangle _{\Lambda_-^2}
     \leq \sum\limits_{i=1}^3 \lambda_i^+ b_i^2 + 
        \sum\limits_{i=1}^3 \lambda_i^- b_i^2 .
\tag4.23
$$
Combining (4.22) and (4.23) we obtain (4.21).
\hskip 4.3cm $\square$

On the set $M_+^4$, $g_\infty$ satisfies
$$
\sigma_2(A) = {\tsize{1\over 4}} | W|^2 = \Vert W \Vert^2 
  = \Vert W^+\Vert^2 + \vert W^- \Vert ^2 .
$$
Therefore,
$$
R = \sqrt{12 |E|^2 + 24 \Vert W^+\Vert ^2
     + 24\Vert W^-\Vert ^2 }.
\tag4.24
$$
Combining (4.21) and (4.24), the integrand
in (4.9) at the point $P$ satisfies the inequality
$$
\align
&-108 \det W^+ - 108 \det W^- + 6 tr E^3
   -9 W_{ij\kappa \ell} E_{i\kappa} E_{j \ell} +
      R\bigl(|E|^2 + {\tsize{3\over 4}} |W|^2 \bigr) \\
& \geq - 108 \lambda_1^+ \lambda_2^+ \lambda_3^+ -
   108 \lambda_1^- \lambda_2^- \lambda_3^- -
   144 b_1b_2b_3 \\
&\qquad\qquad- 36 \bigl[\lambda_1^+ b_1^2 + \lambda_2^+ b_2^2 +
       \lambda_3^+ b_3^2  + \lambda_1^- b_1^2 + \lambda_2^- b_2^2
          + \lambda_3^- b_3^2 \bigr] \\
&+ \Bigl\{ 48\bigl(b_1^2 + b_2^2 + b_3^2 \bigr)  
      + 24\Bigl[ \bigl(\lambda_1^+ \bigr)^2 +
       \bigl(\lambda_2^+\bigr)^2 + \bigl(\lambda_3^+\bigr)^2
   +\bigl(\lambda_1^-\bigr)^2 + \bigl(\lambda_2^-\bigr)^2 +
         \bigl(\lambda_3^- \bigr)^2 \Bigr]\Bigr\}^{{1\over 2}} \\
&\qquad \times \Bigl\{ 4\bigl(b_1^2+b_2^2+b_3^2 \bigr) + 3 
    \Bigl[ \bigl(\lambda_1^+\bigr)^2
    +\bigl(\lambda_2^+\bigr)^2 + \bigl(\lambda_3^+\bigr)^2 +
      \bigl(\lambda_1^- \big)^2 + \bigl(\lambda_2^- \bigr)^2
        +\bigl(\lambda_3^- \bigr)^2 \Bigr]\Bigr\} \\
&\equiv F\bigl(\lambda_1^+, \lambda_2^+ ,\lambda_3^+,
     \lambda_1^-, \lambda_2^- ,\lambda_3^-, b_1,b_2,b_3 \bigr) \,.
\endalign
$$
\enddemo

\proclaim{Proposition 4.4}
Suppose $0 \leq b_1 \leq b_2 \leq b_3$,
$\lambda_1^\pm \leq \lambda_2^\pm \leq \lambda_3^\pm$
with $\si3 \lambda_i^\pm =0$ and 
$\si3 \bigl(|\lambda_i^+ |^2 + |\lambda_i^-|^2 \bigr) \neq 0$.
Then $F(\lambda_1^+, \lambda_2^+, \lambda_3^+,$
$\lambda_1^-, \lambda_2^-, \lambda_3^-, b_1, b_2, b_3 ) \geq 0$, and
equality holds if and only if one of the following is true:
\roster
\item $b_1=b_2 =b_3=0$ and there exists some $a\geq 0$ with
$\lambda_1^+ = \lambda_2^+ = - a, \;\lambda_3^+ = 2a$,
$\lambda_1^-= \lambda_2^- =\lambda_3^- =0$; or
$\lambda_1^+ = -2a, \; \lambda_2^+ = \lambda_3^+ =a$,
$\lambda_1^- = \lambda_2^- = \lambda_3^- =0$; or similar
cases with the role of $\lambda_i^+$ and $\lambda_i^-$ interchanged.
\item $b_1= b_2 = b_3$, $\lambda_i^\pm =0$ for all $1 \leq i \leq 3$.
\endroster
\endproclaim

The proof of Proposition 4.4 is given in the Appendix, and amounts to a complicated 
Lagrange-multiplier problem.  We will assume the result for now, and explain 
how Proposition 4.1
follows.

By Proposition 4.4, the integrand in (4.9) is non-negative.  Since the integral
 is less than or
equal to zero, it follows that the integrand $F(\lambda_1^+, \lambda_2^+, \lambda_3^+,
\lambda_1^-, \lambda_2^-, \lambda_3^-, b_1, b_2, b_3 ) \equiv 0$.  
 Thus, at each point in
$M^4_{+}$ either case (1) or case (2) of Proposition 4.4  above must hold.  Since by definition $|W| > 0$ on $M^4_{+}$,
case (1) is the only possibility.  In particular, $E \equiv 0$ on $M^4_{+}$ and at each point 
either $W^{+} =0$ or $W^{-} = 0$.

Since $E \equiv 0$ on $M^4_{+}$ the scalar curvature
is constant on each component of $M^4_{+}$, which implies by (4.9) that $|W|^2$ is also constant on 
each component. 
We claim that $M^4_{+} = M^4$; i.e., $M^4_{0}$ is empty.  To see why, choose a component
$O$ of $M^4_{+}$ and a sequence of points $\{x_i\}$ in $O$ with $x_i \to x_0 \in M^4_{0}$.
Since $|W|^2$ is constant in $O$, 
$$
c = |W_{g_{\infty}}|^2 (x_i)
$$
for some $c>0$.  By conformal invariance of the Weyl tensor, 
$$
c =  |W_{g_{\infty}}|^2(x_i) = e^{-4w_{\infty}(x_i)}|W_{g}|^2 (x_i).
$$
By definition, $|W_{g}|^2 (x_0) = 0$, and consequently $w_{\infty}(x_i) \to -\infty$ as $i \to \infty$.
But this contradicts the fact that $w_{\infty} \in C^{1,1}$.  It follows that the Weyl tensor
cannot vanish on $M^4$, so $(M^4,g_{\infty})$ is a smooth Einstein manifold.  Moreover, since $|W|^2$ is constant
and either $W^{+} =0$ or $W^{-} = 0$ at each point, it follows that one of the components of the Weyl 
tensor vanishes identically on $M^4$.
By Hitchin's classification result [Hi], $(M^4,g_{\infty})$
is homothetically isometric to  $\pm \bold C P^2$ with the Fubini-Study metric.  This completes the
proof of Proposition 4.1.
 

\head {\bf Appendix }\endhead


In this appendix we establish Proposition A below, which
is slightly more general than Proposition 4.4.

Denote $\vec{B}= (b_1,b_2,b_3)$, $\vec{X} = (x_1, x_2, x_3)$,
$\vec{Y} = (y_1, y_2, y_3)$ vectors in $\Bbb R^3$, and
$| \vec{B}|^2 = \sum_{i=1}^3 b_i^2$, $| \vec{X}|^2 = \si3 x_i^2$,
$|\vec{Y}|^2 = \si3 y_i^2$.  Define the functional
$I=I\bigl(\vec{B}, \vec{X}, \vec{Y}\bigr)$ by
$$
\align
I&= \sqrt{6} \; \bigl[4 \vba + 3\bigl(\vxa +\vya \bigr)\bigr] 
    \; \bigl(2 \vba + \vxa + \vya \bigr)^{1/2} \\
& \quad - 54 x_1x_2x_3 - 54 y_1y_2y_3 - 72 b_1b_2b_3 \\
&\quad -18\bigl(x_1b_1^2 + x_2b_2^2 + x_3 b_3^2 + y_1 b_1^2 
    +y_2 b_2^2 +y_3b_3^2 \bigr)\;.
\endalign
$$
\proclaim{Proposition A}
Assume $b_1^2 \leq b_2^2 \leq b_3^2 , \si3 x_i =0 , \si3 y_i = 0$.
Then $I \geq 0$.  Further, $I=0$ only at the following points: 
\roster
\item"{(i)}"  $\vec{B} = \vec{X} = \vec{Y} = (0,0,0), \quad \text{ or}$
\item"{(ii)}"  $\vec{B} = (0,0,0)$ and either $\vec{X} = (-a,-a,2a)$,
$\vec{Y} = (0,0,0)$, or a permutation of $x_1 = -a, x_2 = -a, x_3 = 2a$ and
$\vec{Y} = (0,0,0)$, or with the preceding values with the roles of $\vec{X}$ and $\vec{Y}$ reversed 
for some $a \neq 0$, or
\item"{(iii)}" $\vec{B} = (b,b,b)$ for some $b\neq 0$ and
$\vec{X} = \vec{Y} = (0,0,0)$. 
\endroster
\endproclaim
\remark{Remark}
If we set $F = 2I$, $x_i = \lambda_i^+$, $ y_i = \lambda_i^-$ 
($1 \leq i \leq 3$), with
$\vxa + \vya \neq 0$, then Proposition 4.4 is a consequence of Proposition A.

We will first establish
Proposition A in the special case where $\vec{B} = (0,0,0)$.
\endremark
\proclaim{Lemma 1}
Denote 
$$
\align
&J\bigl(\vec{X} , \vec{Y} \bigr) =I \bigl(0, \vec{X} , \vec{Y} \bigr)\\
& \qquad = 3 \sqrt{6} \bigl( \vxa + \vya ) (\vxa + \vya \bigr)^{1/2} \\
&\qquad\qquad -54 x_1x_2x_3 - 54 y_1y_2y_3
\endalign
$$
with $\si3 x_i = \si3 y_i = 0$. Then $J \geq 0$, and $ J =0$ only when
$\vec{X} = \vec{Y} = (0,0,0)$ or at the points
$$
\align
\vec{X} &= (-a, -a, 2a) , \quad \vec{Y} = (0,0,0), \; \text{or} \\
\vec{X} &= (2a, -a, -a) , \quad \vec{Y} = (0,0,0), \; \text{or}\\
\vec{X} &= (-a, 2a, -a) , \quad \vec{Y} = (0,0,0); \\
\endalign
$$
or with the roles of $\vec{X}$ and $\vec{Y}$ reversed, for some $a \neq 0$.
\endproclaim
\demo{Proof}
In the case $\vec{X} = \vec{Y} = (0,0,0), \; J=0$.
In the case $\vxa + \vya \neq 0$, replacing $x_i$ and $y_i$ by
${\pm x_i\over \bigl(|x|^2 + |y|^2 \bigr)^{1/2}}$ and 
${\pm y_i\over \bigl(|x|^2 + |y|^2 \bigr)^{1/2}}$,
we may assume w.l.o.g. that
$\vxa +\vya = 1$, $x_1 \leq x_2 \leq x_3$, and 
$y_1 \leq y_2 \leq y_3$.  

To find the minimal points of $J$ we use the method of Lagrange multipliers.
Let $\varphi_1 = \vxa +\vya = 1$,  $\varphi_2 = x_1 + x_2 + x_3 =0$, and
$\varphi_3 = y_1 +y_2 + y_3 =0$ denote the constraints, and let
$\mu , 2 \beta, 2\gamma$ denote the respective Lagrange multipliers. 
Set 
${\partial J\over \partial x_i}=\mu {\partial \varphi_1\over \partial x_i}$
$+ 2 \beta {\partial \varphi_2\over \partial x_i}$,
${\partial J\over \partial y_i} = \mu{\partial \varphi_1\over \partial y_i}$
$+ 2\gamma  {\partial \varphi_3\over \partial y_i}$
for $i=1,2,3$; we get the equations
\TagsOnLeft
$$
\align
-27 \; x_2x_3 &= \mu\; x_1 + \beta  
\tag{\text{A1}}\\
-27 \; x_1 x_3 &= \mu \;x_2 + \beta 
\tag{\text{A2}} \\
-27 \; x_1 x_2 &= \mu \;x_3 + \beta 
\tag{\text{A3}} \\
\endalign
$$
Subtracting (A2) from (A1) we get
$$
-27(x_2 -x_1)x_3 = \mu (x_1 - x_2),
$$
so we have either $x_1 = x_2$, or $x_1 \neq x_2$ while
$27 x_3 = \mu$.  
Subtracting (A3) from (A2) we get similarly $x_2 = x_3$ or
$x_2 \neq x_3$ while $27x_1 = \mu$.
Thus, we have three possibilities: either
$x_1 = x_2 = x_3 = 0$, or $x_1 = x_2 \neq x_3$ with
$27x_3 = \mu$, or $x_1 \neq x_2 = x_3$ with
$27 x_1 = \mu$.
In summary, we have $\vec{X} = (0,0,0)$ or $\vec{X} = (-a,-a,2a)$
with $54a = \mu$; or $\vec{X} = (2a,-a,-a)$ with $54a= \mu$,
where $a \geq 0$.  By symmetry, we also have either $\vec{Y} = (0,0,0)$,
or $\vec{Y} = (-c,-c,2c)$ with $54c= \mu$, or
$\vec{Y} = (-2c, c, c)$ with $-54c = \mu $, where $c \geq 0$.

Combining the possibilities for $\vec{X}$ and $\vec{Y}$,
we have eight cases.
\roster
\item"{(i)}"  $\vec{X}=(0,0,0), \quad \vec{Y}=(-c,-c,2c)$; \quad then $J=0$.
\item"{(ii)}"  $\vec{X}=(0,0,0), \quad \vec{Y}=(-2c,c,c)$; \quad  then $J=0$.
\item"{(iii)}"  $\vec{X}=(-a,-a,2a), \quad \vec{Y}=(-c,-c,2c)$ with
$54a= \mu = 54 c \geq 0$; \quad then $J (\vec{X} , \vec{Y} ) =$
$3\sqrt{6} \cdot 12a^2$ $\sqrt{12} \, a - 54 \cdot 4a^3 =$
$216\bigl(\sqrt{2} -1 \bigr) a^3 \geq 0$. 
\item"{(iv)}"  $\vec{X} = (-a, -a, 2a)$, $\vec{Y} = (-2c, c, c)$ with
$54a= \mu = -54c$;  \quad then $a =-c$ while both $a \geq 0, c\geq 0$;
thus $\vec{X} = \vec{Y} = (0,0,0)$ and $J=0$.
\endroster
Cases (v) $\rightarrow$ (viii) are similar to cases (ii) to (vi)
with the roles of $\vec{X}$ and $\vec{Y}$ reversed.  This finishes the 
proof of Lemma 1.
\qed
\enddemo

We now consider the general case in Proposition A.
The proof is more tedious, but follows the same pattern of the
proof of Lemma 1.  We first outline the steps.

\demo{Outline of the proof of Proposition A.}
When $\vec{B} = (0,0,0)$, we apply Lemma 1.
So assume $\vec{B} \neq (0,0,0)$.  We also assume w.l.o.g. that
$ 2 \vba + \vxa + \vya = 1$.  To locate the minimal points of $I$
under this constraint, we once again apply the method of Lagrange multipliers.
We will break the proof into the following four steps:

\noindent $\underline{\text{Step 1}}$.  We may assume w.l.o.g. that
$0 \leq b_1 \leq b_2 \leq b_3$ and
$x_1 \leq x_2 \leq x_3$, $y_1 \leq y_2 \leq y_3$, and
$x_3 \geq 0 , y_3 \geq 0$.

\noindent $\underline{\text{Step 2}}$.  Actually, $b_1 = b_2 = b$.

\noindent $\underline{\text{Step 3}}.$
When $b=0$, then $x_1 = x_2 = y_1 = y_2$ and $I >0$.

\noindent $\underline{\text{Step 4}}.$
When $b >0$, then $x_1 = x_2 = y_1 = y_2 = -a$,
and either $a=0$, $b_3 = b$, and $I= 0$; \quad or $a>0$
and $I > 0$.

We now prove each of the steps in more detail.
\enddemo

\demo{Proof of Step 1}
We first observe that at a minimal point of 
$I=I\bigl(\vec{B}, \vec{X}, \vec{Y}\bigr)$
we may assume $b_1b_2b_3 \geq 0$.  Thus by switching
$b_i$ with $-b_i$ $(i=1, 2,3)$, we may assume under the hypothesis
$b_1^2 \leq b_2^2 \leq b_3^2$ that
$0 \leq b_1 \leq b_1 \leq b_3$.

Next observe for any $\vec{B}$, for $\vec{X}$ to be a minimal point of
$I$ we must have $x_1b_1^2 + x_2b_2^2 + x_3b_3^2 \leq$
$x_2b_1^2 + x_1b_2^2 +x_3b_3^2$.  Therefore, 
$(x_1-x_2)b_1^2 \leq (x_1-x_2) b_2^2$.  Thus unless
$b_1^2 = b_2^2 =0$, we have
$x_1 \leq x_2$; but when $b_1^2 = b_2^2 =0$, we may also
switch the order of $x_1$ and $x_2$ if necessary, and assume $x_1 \leq x_2$
to attain the same value of $I$.  Thus we can argue
similarly and obtain $x_1 \leq x_2 \leq x_3$ and
$y_1 \leq y_2 \leq y_3$.  Since $\si3 x_i = \si3 y_i =0$, it follows
that $x_3 \geq 0 , y_3 \geq 0$.
\enddemo
\demo{Proof of Step 2}
We now set up the Lagrange multiplier problem under the
constraints
$$
\align
\varphi_1 = 2\vba + \vxa + \vya& = 1\\
\varphi_2 = x_1 + x_2 + x_3 &= 0\\
\varphi_3 = y_1 + y_2 + y_3 &= 0 \\
\endalign
$$
with multipliers $\mu, 2 \beta , 2 \gamma$, respectively.
To locate the minimal point(s) of $I$ under the given constraints
we set
${\partial I\over \partial b_i}=\mu {\partial \varphi_1\over \partial b_i}$,
${\partial I\over \partial x_i}= \mu{\partial \varphi_1\over \partial x_i}$
$+2 \beta {\partial \varphi_2\over \partial x_i}$,
${\partial I\over \partial y_i} = \mu{\partial \varphi_1\over \partial y_i}$
$+2 \gamma {\partial \varphi_3\over \partial y_i}$
for $i =1, 2,3$ and obtain the following nine equations:
$$
\xalignat 2
&-9(x_1+y_1)\; b_1-18 b_2b_3 =(\mu -2 \sqrt{6})b_1 \tag{\text{A4}}\\
&-9(x_2+y_2) \; b_2-18 b_1b_3 =(\mu -2 \sqrt{6})b_2 \tag{\text{A5}}\\
&-9(x_3+y_3) \; b_3-18 b_1b_2 =(\mu -2 \sqrt{6})b_3\tag{\text{A6}}\\
& \qquad -27(x_2x_3) -9 b_1^2 =(\mu -3\sqrt{6})x_1+\beta\tag{\text{A7}}\\
& \qquad -27(x_1x_3)-9 b_2^2 =(\mu -3 \sqrt{6})x_2+\beta\tag{\text{A8}}\\
& \qquad -27(x_1x_2)-9 b_3^2 =(\mu -3\sqrt{6})x_3 + \beta\tag{\text{A9}}\\
\endxalignat
$$
and (A10), (A11), and (A12) are obtained by substituting $y_i$ in place of $x_i$ $(i=1, 2, 3)$ 
and $\gamma$ in place of $\beta$
in equations
(A7), (A8), (A9).

We now assume $0 \leq b_1 \leq b_2 \leq b_3$,
$x_1 \leq x_2 \leq x_3$, $y_1 \leq y_2 \leq y_3$, and prove that
(A4)---(A12) imply $b_1 = b_2$.
To see this, first suppose $b_1 =0$.  Then by (A4),
$b_2 b_3 =0$; since $b_2 \leq b_3$, $b_2 =0$ and thus $b_1 = b_2$.

If $b_1 \neq 0$, then $b_2 \neq 0$, $b_3 \neq 0$.
Subtracting (A8) from (A7) we get
\TagsOnLeft
$$
 -27 x_3 (x_2-x_1) - 9(b_1^2 -b_2^2) =
    (\mu -3 \sqrt{6}) (x_1 - x_2) .
\tag{\text{A13}}
$$
Subtracting (A11) from (A10) we get
$$
-27 y_3 (y_2-y_1) - 9(b_1^2 -b_2^2)
  = (\mu - 3 \sqrt{6}) (y_1 -y_2).
\tag{\text{A14}}
$$
Adding (A13) and (A14),
$$
\align
-18 (b_1^2 - b_2^2 ) &=(\mu -3 \sqrt{6}) (x_1+y_1-x_2-y_2)\\
  & \quad +27 \bigl(x_3(x_2-x_1)+y_3 (y_2-y_1)\bigr).
\tag{\text{A15}}
\endalign
$$

Subtracting (A5) from (A4) and substituting into
(A15), we obtain
$$
-18 (b_1^2 - b_2^2) = 2 (\mu - 3\sqrt{6})
   \left(-{b_2b_3\over b_1} + {b_1b_3\over b_2} \right) + 27 P
\tag{\text{A16}}
$$
where $P = x_3 (x_2 - x_1) + y_3 (y_2 - y_1)$.
Observe that $P \geq 0$.

Thus if $b_2^2 - b_1^2 > 0$, we may divide (A16) by
$b_2^2 - b_1^2$ and get
$$
\quad 18 = -2( \mu -3\sqrt{6}) {b_3\over b_1b_2} 
     \; + 27 {P\over b_2^2 - b_1^2}.
\tag{\text{A17}}
$$
Substituting (A17) into (A6), we get
$$
\align
-9(x_3+y_3) &=
   (\mu -3\sqrt{6} ) + \sqrt{6} + 18 {b_1b_2\over b_3} \\
&= 9 {b_1b_2\over b_3 }+ {27\over 2} \; {P\over b_2^2 - b_1^2} \;
    {b_1 b_2\over b_3} + \sqrt{6} .
\tag{\text{A18}}
\endalign
$$
The left--hand side of (A18) is $\leq 0$, while the right--hand side
is $\geq \sqrt{6} >0$.  Since this contradicts
the hypothesis $b_2^2 - b_1^2 > 0$, we must have 
$b_2^2 = b_1^2$ and hence $b_1 = b_2$.
This establishes step 2.
\enddemo 

\demo{Proof of Step 3}
Denote $b_1 = b_2 = b$, assume $b = 0$, and
$b_3 > 0$.

We rewrite (A6),(A13) and (A14) in this case and get
$$
\xalignat 2
-9 (x_3 +y_3 ) & = (\mu - 2 \sqrt{6}), \tag"(\text{A6})$^\prime$"\\
-27 x_3 (x_2 - x_1 ) &= ( \mu - 3\sqrt{6} ) (x_1 - x_2),
     \tag"(A13)$^\prime$"\\
-27 y_3 (y_2-y_1 ) &= (\mu - 3\sqrt{6})(y_1-y_2 ). \tag"(A14)$^\prime$"
\endxalignat
$$
We observe that from (A6)$^\prime$ (and $x_3 \geq 0$, $y_3 \geq 0$)
that $\mu -2 \sqrt{6} \leq 0$.  Thus we conclude from
(A13)$^\prime$ and (A14)$^\prime$ that $x_1 = x_2$ and $y_1 = y_2$.
Denote $x_1 = x_2 = -a$, $y_1 = y_2 = -c$, with $a,c \geq 0$.
Then $\vec{B} = (0,0,b_3)$, $\vec{X} = (-a, -a, 2a)$, 
$\vec{Y} = (-c, -c, 2c)$.  We now assert that $a=c$.  To see this,
we first subtract (A9) from (A8) to get
$$
\align
&\qquad -27 x_1 (x_3 - x_2 ) + 9 b_3^2 = ( \mu - 3 \sqrt{6} )
    (x_2 - x_3 )
\tag{\text{A19}} \\
\hskip -2.5cm \text{or} &\\
&\qquad 27a^2 + 3b_3^2 = -(\mu - 3 \sqrt{6}) a\,.
\tag"(A19)$^\prime$"
\endalign
$$
Similarly we subtract (A12) from (A11) and get
$$
27c^2 + 3b_3^2 = -(\mu - 3\sqrt{6}) c \,.
\tag"(A20)$^\prime$"
$$
Subtracting (A20)$^\prime$ from (A19)$^\prime$, we then have
$$
27(a^2 -c^2) = -(\mu - 3 \sqrt{6})(a-c). 
$$
Thus, either $a=c$, or $a \neq c$ while
$$
-27 (a+c) = \mu - 3\sqrt{6} \,.
\tag{\text{A21}} 
$$
We will now derive a contradiction to see that the latter possibility
$(a \neq c)$ does not happen.  Comparing (A21) with (A6)$^\prime$, we see 
$-18(a+c) = \mu -2 \sqrt{6}$.  Combining (A21) with (A6)$^\prime$ we get
$$
a+c = {\sqrt{6}\over 9} \,.
\tag{\text{A22}}
$$
We now add (A19)$^\prime$ to (A20)$^\prime$ and substitute the relation
$\varphi_1 = 2b_3^2 +6(a^2+c^2 )=1$
and (A21) into the equation to get
$$
3(a^2+c^2) = 9(a+c)^2 -1.
\tag{\text{A23}}
$$
By (A22), this expression equals $ {6\over 9} -1 < 0$, which is a contradiction.
Thus we conclude that $a=c$ and $\vec{B} = (0,0,b_3)$,
$\vec{X} = \vec{Y} = (-a, -a, 2a)$
for some $b_3 > 0$ and $a \geq 0$.  At this moment, we can
check directly that $I= I(\vec{B}, \vec{X}, \vec{Y} ) > 0$ if
$b_3 > 0$.  To be more precise, one can check that
$I = 8(I_1 - I_2)$ with $I_1 = \sqrt{3} (b_3^2 + 9a^2)$
$(b_3^2 + 6a^2)^{1/2}$, $I_2 = 9a(3a^2+b_3^2)$, and
$I_1^2 - I_2^2 >0$.  This establishes step 3.
\enddemo

\demo{Proof of Step 4}
Assume $b_1 = b_2 = b \neq 0$, $ b_3 >0$.
In this case, we may write (A4), (A5), (A6) as
$$
\align
-9 (x_1+y_1 ) - 18b_3 &= \mu -2 \sqrt{6}, 
\tag"(A4)$^{\prime\prime}$"\\
-9(x_2+y_2) - 18b_3 &= \mu -2\sqrt{6},
\tag"(A5)$^{\prime\prime}$" \\
-9 (x_3 + y_2) b_3 - 18b^2 & = (\mu - 2\sqrt{6}) b_3 \, .
\tag"(A6)$^{\prime\prime}$"
\endalign
$$
Subtracting (A5)$^{\prime\prime}$ from (A4)$^{\prime\prime}$, we get
$$
x_1+y_1 = x_2 +y_2 \,.
\tag{A24}
$$
Subtracting (A8) from (A7) and (A11) from (A10) we get
$$
\align
-27x_3 (x_2 -x_1) &= (\mu -3 \sqrt{6}) (x_1 - x_2),
\tag"(A13)$^{\prime\prime}$"\\
-27 y_3 (y_2 -y_1 ) & = (\mu - 3 \sqrt{6} ) (y_1 -y_2 )\,.
\tag"(A14)$^{\prime\prime}$"
\endalign
$$
\enddemo

We now make two claims.

\noindent $\underline{\text{Claim 1}}$:  $x_1 = x_2 (= -a)$,
$y_1 = y_2 (= -c)$. 

\noindent $\underline{\text{Claim 2}}$:  $a = c$.

\demo{Proof of Claim 1}
By (A13)$^{\prime\prime}$, we have either $x_1 = x_2$  or 
$x_1 \neq x_2$ and $27 x_3 = \mu - 3 \sqrt{6}$.
By (A14)$^{\prime\prime}$, we have either $y_1 = y_2$
or $y_1 \neq y_2$ and $ 27 y_3 = \mu - 3 \sqrt{6}$.
By (A24), we have $x_1 = x_2$ implies $y_1 = y_2$.
Thus, we either have $x_1 = x_2$ and $y_1 = y_2$ as claimed, or
$x_1 \neq x_2$, $ y_1 \neq y_2$ while $x_3=y_3={1\over 27} (\mu -3\sqrt{6})$.
But in the latter case, from (A6)$^{\prime\prime}$ we would have
$$
-9(x_3 + y_3) b_3 - 18b^2 = 
    (\mu -2 \sqrt{6} )b_3 = (27x_3 + \sqrt{6} )b_3 \,.
$$
This is a contradiction as the left--hand side  of the equation is
less than zero, while the right--hand side is bigger than zero.
\enddemo

\demo{Proof of Claim 2}
We follow the same strategy as in the proof of Step 3.
Subtracting (A9) from (A8) we get
$$
27a^2 + 3(b_3^2 - b^2 ) = -a(\mu - 3\sqrt{6} )\,.
\tag"(A19)$^{\prime\prime}$"
$$
Similarly, if we subtract (A12) from (A11) we get
$$
27 c^2 + 3 ( b_3^2 - b^2) = -c ( \mu - 3\sqrt{6} ) \,.
\tag"(A20)$^{\prime\prime}$"
$$
Finally, subtracting (A20)$^{\prime\prime}$ from (A19)$^{\prime\prime}$ we have
$$
27(a^2 -c^2 ) = - (\mu - 3 \sqrt{6} ) (a-c) \,.
$$
Thus, either $a=c$ as claimed, or
$$
a\neq c \quad \text{ and } \quad 27(a+c) = - (\mu - 3\sqrt{6}) \,.
\tag"(A21)$^{\prime\prime}$"
$$
We will now show that (A21)$^{\prime\prime}$ cannot be true.
To see this, denote $a+c = \ell \, (\ell \geq 0)$, $a-c =d $.
Then $d \neq 0$, and w.l.o.g. we may assume $d > 0$.  
Substituting (A21)$^{\prime\prime}$ into (A4) we get
$$
\align
-18b_3& = (\mu - 2 \sqrt{6} ) - 9 (a+c) \\
   &= (\mu - 3 \sqrt{6} ) - 9 \ell + \sqrt{6} = -36\ell +\sqrt{6} \,.
\tag{A24}
\endalign
$$
Substituting (A21)$^{\prime\prime}$ to (A6) we get
$$
\align
-18b^2 &=  (\mu -2 \sqrt{6} ) b_3 + 18 (a+c) b_3\\
&= (-9\ell + \sqrt{6} ) b_3 \,.
\tag{A25}
\endalign
$$
Thus in particular
$$
-9 \ell + \sqrt{6} < 0 \,.
\tag"(A25)$^{\prime\prime}$"
$$
Combining (A24) and (A25), we get
$$
18(b^2 -b_3^2 ) = (-36\ell + \sqrt{6} + 9 \ell - \sqrt{6} ) b_3 =
   -27 \ell b_3\,.
\tag{A26}
$$
On the other hand, substituting (A21)$^{\prime\prime}$ into (A19)$^{\prime\prime}$,
we get
$$
\align
27a^2 +3 (b_3 - b^2 ) &= 27a(a+c), \\
&\hskip -2.66in \text{so} \\
3(b_3^2 -b^2)& = 27 ac= 27 {\ell^2 -d ^2\over 4} \,.
\tag{A27}
\endalign
$$
Combining (A24), (A26), (A27), we get
$$
27d^2 = ( - 9 \ell + \sqrt{6}) \ell,
\tag{A28}
$$
which contradicts (A25)$^{\prime\prime}$.
We conclude $a=c$, as in Claim 2.

We are now in the situation where $\vec{B} = (b,b,b_3)$, with
$b_3 \geq b > 0$, and $\vec{X} = \vec{Y} = (-a,-a,2a)$.
There are two final possibilities to consider, depending on the 
sign of $a$.

\noindent $\underline{\text{Claim 3}}$ If $a=0$, then
$b_3 = b\neq 0$, and $I \equiv 0$.

\noindent $\underline{\text{Claim 4}}$ If $a>0$ then $I>0$.
\enddemo

\demo{Proof of Claim 3}
When $a=0$, we multiply (A4)$^{\prime\prime}$ by $b_3$ then subtract
(A6)$^{\prime\prime}$ to get $b_3^2 = b^2$, hence $b_3 = b$.
In this case $I= 72b^3 - 72b^3 \equiv 0$.
\enddemo

\demo{Proof of Claim 4}
When $a \neq 0$, we will show that $I > 0$.
First, we rewrite (A4)$^{\prime\prime}$, (A6)$^{\prime\prime}$, 
(A7) and (A9) as follows:
$$
\align
18a - 18b_3 & = \mu - 2\sqrt{6},
\tag"(A4)$^{\prime\prime}$"\\
-36ab_3 - 18b^2 &= (\mu - 2 \sqrt{6} ) b_3,
\tag"(A6)$^{\prime\prime}$"\\
54a^2 -9b^2 &= - (\mu - 3\sqrt{6}) a + \beta,
\tag"(A7)$^{\prime\prime}$"\\
-27a^2 - 9b_3^2 &= 2 (\mu - 3 \sqrt{6}) a+ \beta\, .
\tag"(A9)$^{\prime\prime}$"
\endalign
$$
Multiplying (A4)$^{\prime\prime}$ by $b_3$ and subtracting  
(A6)$^{\prime\prime}$ from the result we get
$$
b_3^2 - b^2 = 3 ab_3\,.
\tag{A29}
$$
Subtracting (A9)$^{\prime\prime}$ from (A7)$^{\prime\prime}$ we get
$$
81a^2 + 9(b_3^2 - b^2 ) = -3 a(\mu - 3\sqrt{6} )\,.
\tag{A30}
$$
Combining (A29) and (A30), we get (for $a \neq 0$)
$$
27a + 9b_3 = - (\mu - 3 \sqrt{6} ) \,.
\tag{A31}
$$
Combining (A4)$^{\prime\prime}$ and (A31) we find
$$
45a - 9b_3 = \sqrt{6} \,.
\tag{A32}
$$
Substituting (A29) into the constraint  $ \varphi_1 = 2(2b^2+b_3^2)+12a^2=1$
we get
$$
6b_3^2 - 12ab_3 + 12a^2 = 1\,.
\tag{A33}
$$
We now introduce the notation $s = a-b_3$ and
rewrite (A32) and (A33) as
$$
\align
36a + 9 s = \sqrt{6},
\tag"(A32)$^\prime$"\\
6(s^2 + a^2) = 1\;.
\tag"(A33)$^\prime$"
\endalign
$$
Applying (A29) we can write $I$ as
$$
\align
I &= \bigl(4(2b^2+b_3^2) + 36a^2\bigr)
   \sqrt{6} - 216a^3 - 72b^2b_3 - 72a(b_3^2 -b^2 )\\
&= 12 \sqrt{6} (b_3^2 -2ab_3 + 3a^2) -72(3a^3 +b_3^2 -3a^2b_3 +3ab_3^2)\\
&= 72 \left[{s^2+2a^2\over \sqrt{6} } - (4a^3 - s^3 ) \right] \,.
\endalign
$$
It remains to check numerically that a solution $(a,s)$ with $a > 0$ of
equations (A32)$^\prime$ and (A33)$^\prime$ is given by 
$a = 0.1617...$, $s =-0.3746...$, and $I = 72(0.079... - 0.069...) > 0$.

We have thus finished the proof of Step 4 and completed the proof of 
Proposition A.
\enddemo



\end{document}